\documentclass[a4paper,10pt,fleqn]{amsart}

\usepackage[utf8]{inputenc}
\usepackage{amsmath}
\usepackage{amsfonts}
\usepackage{amsthm}
\usepackage{amssymb}
\usepackage{amsaddr}
\usepackage{bbm}
\usepackage{xcolor}
\usepackage{hyperref}
\usepackage{nameref}
\usepackage{mathtools}

\title{A framework for non-local, non-linear initial value problems}

\author[G.~Karch, M.~Kassmann, M.~Krupski]{%
Grzegorz Karch\MakeLowercase{$^a$},
Moritz Kassmann\MakeLowercase{$^b$},
Miłosz Krupski\MakeLowercase{$^{ac}$}}
\email{\vtop{\hsize=.707\textwidth\noindent grzegorz.karch@uwr.edu.pl,
moritz.kassmann@uni-bielefeld.de,
milosz.krupski@uwr.edu.pl}}

\address{
 $^a$Instytut Matematyczny, Uniwersytet Wrocławski,\\
 pl. Grunwaldzki 2/4, 50-384 Wrocław, Poland\bigskip\\
 $^b$Universität Bielefeld, Fakultät für Mathematik,\\
 Postfach 10 01 31, 33501 Bielefeld, Germany\bigskip\\
 $^c$Institutt for matematiske fag, NTNU,\\
 7491 Trondheim, Norway}

\thanks{Financial support by the Deutsche Forschungsgemeinschaft (DFG) through CRC 1283 is gratefully acknowledged. The third author was supported by the Polish NCN grant 2016/23/B/ST1/00434.}

\date{\today}

\subjclass[2010]{35K55, 35K65, 35R11, 35L65}
\keywords{Non-linear diffusion, non-local evolution equations}

\newtheorem{theorem}{Theorem}[section]
\newtheorem{lemma}[theorem]{Lemma}
\newtheorem{corollary}[theorem]{Corollary}
\newtheorem{proposition}[theorem]{Proposition}
\theoremstyle{remark}
\newtheorem{remark}[theorem]{Remark}
\theoremstyle{definition}
\newtheorem{definition}[theorem]{Definition}

\DeclareMathOperator{\sgn}{sgn}
\DeclareMathOperator{\loc}{loc}
\let\div\relax
\DeclareMathOperator{\div}{div}
\DeclareMathOperator{\supp}{supp}
\newcommand{\dt} {\partial_t}

\newcommand{\T} {[0,\infty)}
\newcommand{\R} {\mathbb{R}}
\newcommand{\X} {{\R^N}}
\newcommand{\XX} {{\R^{2N}}}
\newcommand{\meas} {\rho}
\renewcommand{\epsilon} {\varepsilon}
\newcommand{\mease} {\meas^{\epsilon}}
\newcommand{\Lap}{\mathcal{L}}
\newcommand{\Lpp}[1]{L^{[1,\infty]}(#1)}
\newcommand{\kernel}[1]{\meas\big(#1(x),#1(y);x,y\big)}
\newcommand{\kernelt}[1]{\meas\big(#1(t,x),#1(t,y);x,y\big)}
\newcommand{\kernell}[2]{\meas(#1,#2;x,y)}
\newcommand{\kernelel}[2]{\mease(#1,#2;x,y)}
\newcommand{\sg}[1]{\eta{(#1)}}

\makeatletter
\let\orgdescriptionlabel\descriptionlabel
\renewcommand*{\descriptionlabel}[1]{%
  \let\orglabel\label
  \let\label\@gobble
  \phantomsection
  \edef\@currentlabel{#1}%
  \let\label\orglabel
  \orgdescriptionlabel{#1}%
}
\makeatother

\makeatletter
\newcommand{\subalign}[1]{%
  \vcenter{%
    \Let@ \restore@math@cr \default@tag
    \baselineskip\fontdimen10 \scriptfont\tw@
    \advance\baselineskip\fontdimen12 \scriptfont\tw@
    \lineskip\thr@@\fontdimen8 \scriptfont\thr@@
    \lineskiplimit\lineskip
    \ialign{\hfil$\m@th\scriptstyle##$&$\m@th\scriptstyle{}##$\crcr
      #1\crcr
    }%
  }
}
\makeatother

\begin{document}
\begin{abstract}
 We study the Cauchy problem with non-linear non-local operators that may 
be degenerate. Our general framework includes cases where the 
jump intensity is allowed to depend on the values of the solution itself,
e.g. the porous medium 
equation with the fractional Laplacian  and the parabolic fractional 
$p$-Laplacian. We show the existence, uniqueness  of 
bounded solutions and study their further properties. Several new examples 
of non-local, non-linear operators are provided.
\end{abstract}

\maketitle
\begin{section}{Introduction}\label{sec:introduction}
 One of the core ideas in describing many phenomena in natural sciences
  is the notion of linear and non-linear 
  diffusion. We refer to  the review article~\cite{MR3588125}
for a~comprehensive introduction to this concept from an analytical perspective.
  The most prominent diffusion is the Brownian Motion, whose mathematical description goes back to the beginning of the 20th century.
  Its infinitesimal generator is the Laplace operator $-\Delta$.

The more general 
class of L\'evy stochastic processes has been found to be important for the modeling of 
diffusive and non-diffusive phenomena~\cite{MR1881970,Humphries2010, schoutens2003levy,MR1833700}. Here, the so 
called rotationally invariant $\alpha$-stable jump process can be seen as an 
important non-local counterpart of the Brownian Motion. Its infinitesimal 
generator is the fractional Laplace operator $(-\Delta)^{\alpha/2}$, where 
$\alpha \in (0,2)$. 
An introduction and overview of some analytical results related to non-local phenomena
can be found in~\cite{MR3469920}.

  In this paper we put forward a framework for non-local, non-linear ``diffusion''.
  The operators that we consider might be degenerate but still 
allow for the conservation of
  mass, the maximum principle, and the comparison principle. It encompasses some of the known
  examples of equations which are used to describe behaviour of non-local diffusive-type processes, like the linear heat 
  equation with the fractional Laplacian, the fractional porous medium equation
  or the parabolic equation with the fractional $p$-Laplacian.
  The framework also allows for new examples to be 
  constructed and studied. 
  
  \subsection*{Statement of the problem}
  Consider the initial value problem
  \begin{equation}\label{pns}
    \left\{
    \begin{aligned}
      &\dt u + \Lap_u u = 0\quad&\text{on $\T\times\X$}, \\
      &u(0) = u_0\quad&\text{on $\X$},
    \end{aligned}
    \right.
  \end{equation}
  where the non-linear non-local operator $\Lap$ is defined by the formula
  \begin{equation}\label{operator}
   \big(\Lap_vu\big)(x) = \int_\X \big[u(x)-u(y)\big] \kernel{v}\,dy
  \end{equation}
  for a given \emph{homogeneous jump kernel} $\meas$.
  
  \begin{definition}\label{def:jump-kernel}
  We say a function $\meas:(\R\times\R)\times(\X\times\X)\to\R$ is a homogeneous jump kernel if it satisfies conditions \ref{a1:positive}--\ref{a6:lipschitz}.
  For all $a,b,c,d\in\R$ and 
  for almost every $(x,y)\in\X\times\X$ we assume that
  \begin{description}
    \item[(A1)\label{a1:positive}] $\meas$ is a non-negative Borel function;
    \item[(A2)\label{a2:symmetry}] $\meas$ is symmetric, i.e.~$\kernell{a}{b} = \meas(b,a;y,x)$;
    \item[(A3)\label{a3:monotonicity}] $\meas$ is monotone in the following sense:
    \begin{equation*}\hspace{\leftmargin}
      (a-b)\kernell{a}{b} \geq (c-d)\kernell{c}{d}\quad\text{whenever $a\geq c\geq d\geq b$};
    \end{equation*}
     \item[(A4)\label{a4:homogeneity}] $\meas$ is homogeneous 
     \begin{equation*}\hspace{\leftmargin}
     \kernell{a}{b} = \meas\big(a,b;|x-y|\big);
    \end{equation*}
     \item[(A5)\label{a5:1levy}]
     for every $R> 0$ there exists a function $m_R:[0,\infty)\to[0,\infty)$ such that
     $\sup_{-R\leq a,b\leq R}\kernell{a}{b}\leq m_R\big(|x-y|\big)$ and
      \begin{equation*}\hspace{\leftmargin}
        \int_\X\big(1\wedge |y|\big)m_R\big(|y|\big)\,dy = K_R <\infty;
      \end{equation*}
     \item[(A6)\label{a6:lipschitz}] 
     $\meas$ is continuous with respect to the first two variables and it is locally Lipschitz-continuous outside diagonals, i.e.~for every $\epsilon >0$ and every $R>\epsilon$ there exists a constant $C_{R,\epsilon}$
     such that
      \begin{equation*}\hspace{\leftmargin}
        \big|\meas(a,b,x,y) - \meas(c,d,x,y)\big| \leq C_{R,\epsilon} \big(|a-c|+|b-d|\big)\, m_R\big(|x-y|\big)
      \end{equation*}
      for every $a,b,c,d\in[-R, R]$ such that $|a-b|\geq\epsilon$ and $|c-d|\geq\epsilon$.
    \end{description}
  \end{definition}
 
 The name \emph{jump kernel} comes from the probabilistic interpretation
 of the role of the operator $\Lap$ in equation~\eqref{pns}. The function $\meas(a,b,x,y)$ 
describes the density of jumps from $x$ to $y$ within a unit time interval. In 
our setup, this density is moreover allowed to depend on the values $u(x)$ and 
$u(y)$ of the solution in place of parameters $a$ and $b$.
\begin{remark}
  In the integrability condition~\ref{a5:1levy}, one may simply \emph{define} the function $m_R$ as
  $m_R\big(|x-y|\big) = \sup_{-R\leq a,b\leq R}\kernell{a}{b}$. In some situations, however,
  it is more convenient to consider a different majorant (cf. Proposition~\ref{lemma:bielefeld}).
 \end{remark}
 \begin{remark}\label{rem:levy}
   If $\mu$ is a non-negative function and the integral $\int_\X\big(1\wedge |y|\big)\mu(y)\,dy$ is finite,
   then we say that $\mu$ is the density of a Lévy measure with low singularity. This corresponds to our assumption on the function $m_R$ in condition~\ref{a5:1levy} in Definition~\ref{def:jump-kernel}. Note 
that, in case of a general Lévy measure $\mu$, the integral
   $\int_\X\big(1\wedge |y|^2\big)\mu(dy)$
   is assumed to be finite. We do not treat this general case here. 
  \end{remark}
 
 Note that in condition \ref{a4:homogeneity} and in the sequel we allow 
for some ambiguity with respect to the function $\meas$.
 Moreover, for clarity of presentation we abbreviate the notation for the jump kernel and write
  \begin{align*}
  \kernel{u} = \meas_{u,x,y}& &\text{and}& &\kernel{u}\,dy\,dx = d\meas_u.
  \end{align*}
 
   \begin{remark}[Maximum principle]
      An  operator $\big(A, D(A)\big)$  satisfies the positive maximum principle if
  for every  $u \in  D(A)$
  the fact that $u(x_0)=\sup_{x\in\X} u(x)\geq 0$
  for some $x_0\in \X$ 
  implies $A u(x_0) \leq 0$.
  Because of condition~\ref{a1:positive} and formula~\eqref{operator},
  the operator $-\Lap_v$ 
  has this property 
  for each $v\in L^\infty(\X)$ as long as it is supplemented with a suitable domain
  (for example, the $BV$ space, see Lemma~\ref{lemma:bv} below).
  \end{remark}

\subsection*{Main results}
Our goal is to prove the results gathered in Theorem~\ref{thm:main}
and Corollary~\ref{cor:main}.

\begin{theorem}\label{thm:main}
Let $\meas$ be a homogeneous jump kernel in the sense of Definition~\ref{def:jump-kernel}.
For every initial condition $u_0\in L^\infty(\X)\cap BV(\X)$,
problem \eqref{pns} has a unique very weak solution $u$ such that
\begin{equation*}
u \in L^\infty \big([0,\infty), L^\infty(\X)\cap BV(\X) \big)\cap W^{1,1}_{\loc}\big([0,\infty),L^1(\X)\big)
\end{equation*}
(see Definition~\ref{def:bvsolution}). This solution has the following properties
\begin{itemize}
\item mass is conserved: $\int_{\X} u(t,x)\,dx = \int_{\X} u_0(x)\,dx$ for all $t\geq 0$;
\item $L^p$-norms are non-increasing: $\|u(t)\|_p\leq \|u_0\|_p$ for all $p\in [1,\infty]$ and $t\geq 0$;
\item if $u_0(x)\geq 0$ for almost every $x\in\X$
      then $u(t,x)\geq 0$ for almost every $x\in \X$ and $t\geq 0$.
\end{itemize}
Moreover, for two solutions $u$ and $\widetilde u$ corresponding to initial conditions $u_0$ and
$\widetilde u_0$, respectively, we have
\begin{equation*}
 \|u(t)-\widetilde u(t)\|_1\leq \|u_0-\widetilde u_0\|_1\quad\text{for every $t\geq 0$}
 \end{equation*}
 and 
 if $u_0(x)\geq \widetilde u_0(x)$ for almost every $x\in\X$
      then $u(t,x)\geq \widetilde u(x,t)$ for almost every $x\in \X$ and $t\geq 0$.
\end{theorem}
\begin{corollary}\label{cor:main}
 Let $\meas$ be a homogeneous jump kernel in the sense of Definition~\ref{def:jump-kernel}.
 For every initial condition $u_0\in L^1(\X)\cap L^\infty(\X)$,
 problem \eqref{pns} has a very weak solution $u$ such that
 \begin{equation*}
 u \in L^\infty \big([0,\infty), L^1(\X)\cap L^\infty(\X)\big)\cap C\big([0,\infty), L^1(\X)\big).
 \end{equation*}
 This solution has the following properties
\begin{itemize}
\item mass is conserved: $\int_{\X} u(t,x)\,dx = \int_{\X} u_0(x)\,dx$ for all $t\geq 0$;
\item $L^p$-norms are non-increasing: $\|u(t)\|_p\leq \|u_0\|_p$ for all $p\in [1,\infty]$ and $t\geq 0$;
\item if $u_0(x)\geq 0$ for almost every $x\in\X$ then $u(t,x)\geq 0$ for almost every $x\in \X$ and $t\geq 0$.
\end{itemize}  
\end{corollary}
\begin{remark}
 Notice that in the statement of Corollary~\ref{cor:main} we consider more general initial data and consequently
 there is no claim of uniqueness of solutions nor a notion of $L^1$-contraction.
\end{remark}

\begin{remark}
One could expect solutions to diffusion equations to be more regular than their initial data (cf. the heat equation). In general, this 
is not the case in our framework and can easily be disproved by the trivial example  $\kernell{a}{b}= 0$.
It satisfies Definition~\ref{def:jump-kernel}, but there is no smoothing effect in the Cauchy problem
\begin{equation*}
  \partial_t u =0, \qquad u(0)=u_0.
\end{equation*}
This example also shows that we cannot expect any decay of solutions.
\end{remark}
\subsubsection*{Strategy of the proof of Theorem~\ref{thm:main}}
Let us briefly describe the strategy to prove Theorem~\ref{thm:main} and thus also the general outline of the paper.
In Section~\ref{sec:regular}, Definition~\ref{def:regular}, we introduce
the notion of \emph{regular jump kernels} in order to 
 show in Theorem~\ref{thm:classical} that in such cases there
exist unique classical solutions to problem~\eqref{pns} for every initial condition $u_0\in L^1(\X)\cap L^\infty(\X)$.
Then, in Section~\ref{sec:existence}, we ``regularize'' homogeneous jump kernels (Lemma~\ref{lemma:approx}). In Theorem~\ref{thm:existence},
we show that if $u_0\in BV(\X)\cap L^\infty(\X)$ then the approximating sequence of solutions
$u^{\epsilon}$ to problems~\eqref{eq:approx}, with initial conditions $u_0$, corresponding to a sequence of regularized
jump kernels, has a convergent subsequence and its limit $u$ 
is a strong solution to problem~\eqref{pns} (see Definition~\ref{def:bvsolution}).
The solution satisfies
$u\in L^\infty\big([0,\infty),BV(\X)\cap L^\infty(\X)\big)\cap C\big([0,\infty),L^1_{\loc}(\X)\big)$.
In Theorem~\ref{thm:contraction} we prove the $L^1$-contraction property, which
immediately implies uniqueness of solutions (Corollary~\ref{cor:uniqueness}).
The properties of solutions indicated in the statement of Theorem~\ref{thm:main}
follow as side-effects in the process of 
this construction (see Corollaries~\ref{cor:Lp-weak}, \ref{cor:mass_conservation} and~\ref{cor:positivity}).
\subsubsection*{Strategy of the proof of Corollary~\ref{cor:main}}
As a consequence of Theorem~\ref{thm:main}, we obtain Corollary~\ref{cor:main}. In Theorem~\ref{thm:existence2} we consider general initial condition $u_0\in L^1(\X)\cap L^\infty(\X)$ and a sequence of its approximations for which we may
use Theorem~\ref{thm:main} to construct unique strong solutions. We are able to verify Definition~\ref{def:distributional}
of the pointwise limit of these solutions
and to deduce properties of this solution
as a consequence of the pointwise convergence of approximations (see Corollary~\ref{cor:properties}).
\subsection*{Scope of the framework}
Let us introduce our framework by presenting several examples. Many 
of them have been discussed intensively in the literature but some of them, to 
the best of our knowledge, are new.
A more detailed description of these examples follows also in Section~\ref{sec:examples}.

\subsubsection*{Known models}
As a simple example which could be written in the form \eqref{pns} we recall the linear
fractional heat equation 
\begin{equation}\label{eq:heat}
\partial_t u +(-\Delta)^{\alpha/2} u=0 \quad \text{with the jump kernel} \quad
\kernell{a}{b}=\frac{C_{\alpha,N}}{|x-y|^{N+\alpha}}.
\end{equation}
Because of the integrability constraint~\ref{a5:1levy}, in this work we have to
assume $\alpha\in (0,1)$. We refer to~\cite{MR3469920} for a gentle introduction
to the fractional Laplacian. Equation~\eqref{eq:heat} has an explicit solution 
$u(t)=p_\alpha(t)*u_0$ with the $\alpha$-stable density
$p_\alpha(t) \in C^\infty(\X)\cap L^1(\X)\cap L^\infty(\X)$ (see~\cite{MR2373320,MR3211862}).

A more complicated case is the fractional porous medium equation   
\begin{equation}\label{fpme}
\partial_t u +(-\Delta)^{\alpha/2} \big(|u|^{m-1}u\big)=0,
\end{equation}
which is studied in depth in \cite{MR2334594,MR2737788,MR2954615}. See also~\cite{MR3666562} for most recent results
and more references.
In our case, we can write this equation as~\eqref{pns}, using the following jump kernel with
$\alpha\in(0,1)$ and $m\geq1$
\begin{equation*}
 \kernell{a}{b}=\frac{|a|^{m-1}a-|b|^{m-1}b}{a-b}\cdot\frac{C_{\alpha,N}}{|x-y|^{N+\alpha}}.
\end{equation*}

Let $\mu$ be a Lévy measure as defined in Remark~\ref{rem:levy}. One can replace the operator
$(-\Delta)^{\alpha/2}$ in equation~\eqref{fpme} by a more general \emph{Lévy
operator}
\begin{equation}\label{L:lin}
\mathcal{L}^\mu\psi(x) 
= \int_{\X\setminus\{0\}} \Big(\psi(x+y) - \psi(x) - y\cdot\nabla\psi(x)\mathbbm{1}_{|y|\leq 1}(y)\Big)\,d\mu(y),
\end{equation}
and the function $u\mapsto |u|^{m-1}u$ by $u\mapsto f(u)$ which is (for example) a continuous, increasing
function with $f(0)=0$. Some modifications of equation~\eqref{fpme} in this fashion are described 
in~\cite{MR3485132,MR3656476}. The general case is studied in~\cite{MR3724879,MR3570132}. 

Assuming that $f$ is continuous and 
non-decreasing and $\mu$ is a general symmetric L\'{e}vy measure, the authors 
establish a uniqueness result for bounded distributional solutions and several 
estimates. In our framework we are able to prove some of these 
results. 
Our approach is less general than the one in \cite{MR3570132} because 
conditions \ref{a5:1levy} and \ref{a6:lipschitz} limit its scope.

A different but related stream of research focuses on another type of non-linear non-local operator, the 
$s$-fractional $p$-Laplacian (see~\cite{MR3148135,MR3491533,MR3456825} and the references therein) in the context of the following evolution equation
\begin{equation}\label{plap}
\partial_t u +\int_\X \frac{\Phi\big(u(x)-u(y)\big)}{|x-y|^{N+ps}}\,dy =0,\quad \Phi(z) = z|z|^{p-2}.
\end{equation}
Usually it is assumed that $s\in(0,1)$ and $p\in(1,\infty)$. The case $p=2$ reduces to the
linear equation~\eqref{eq:heat}.
For some pairs $(s,p)$, but also other functions $\Phi$, equation~\eqref{plap} can be written in the form of 
equation~\eqref{pns},
see Proposition~\ref{lemma:fpl}.

Our results also apply to a regular version of equation \eqref{plap} i.e.~where
the kernel $|x-y|^{-(N+ps)}\,dy$ is replaced by a sufficiently regular, integrable  and non-negative function $J\big(|x-y|\big)$. This case was studied in 
\cite[Chapter 6]{MR2722295}.

Results  involving regular jump kernels (see Definition~\ref{def:regular} and the entire Section~\ref{sec:regular})
can be directly applied to the following non-local equation studied in~\cite{2016arXiv160203522S}
\begin{equation*}
 \dt u =\int_\X k \big(u(t, x), u(t, y)\big)\big[u(t, y) - u(t, x)\big]J(x - y)\,dy.
\end{equation*}
Here, the kernel $J:\X\to\R$ is an integrable, non-negative function supported in the
unit ball. The function $k : \R^2 \to \R$ is locally
Lipschitz-continuous and non-negative. We refer the reader to the work~\cite{2016arXiv160203522S} for other properties of solutions such as the strong maximum principle.

Models combining local counterparts of operators~\eqref{plap} and~\eqref{fpme}\footnote{We thank Félix del Teso for
signalling this example and pointing us to the references.}, namely
\begin{equation}\label{eq:doubly-nonlinear}
 \dt u = (\Delta_p)\big(u\,|u|^m\big),\quad\text{where}\quad\Delta_pv = \div\big(|\nabla v|^{p-2}\nabla v\big),
\end{equation}
have been studied as well~\cite{MR1097286,MR2988757}. 

\subsubsection*{New models}
Our results can be applied to non-linear, non-local evolution equations which have not been previously studied.
We may combine equations~\eqref{fpme},~\eqref{L:lin} and ~\eqref{plap} and study the following
non-local counterpart of equation~\eqref{eq:doubly-nonlinear}
\begin{equation*}
 \dt u +\int_\X {\Phi\big[f\big(u(x)\big)-f\big(u(y)\big)\big]}\mu\big(|x-y|\big)\,dy = 0.
\end{equation*}
Here $f$ and $\Phi$ are non-decreasing functions and $\mu$ is a density of a Lévy measure with low singularity.
All assumptions are indicated in Proposition~\ref{lemma:doubly-nonlinear}.

By considering what we call \emph{convex diffusion operator} we introduce the following evolution equation
\begin{equation*}
 \dt u + \int_\X \big[u(x)-u(y)\big]\big[f\big((u(x)\big)+f\big(u(y)\big)\big]\mu\big(|x-y|\big)\,dy = 0
\end{equation*}
for a non-negative, convex function $f$ and a density of a Lévy measure with low singularity $\mu$. We discuss this example in
Proposition~\ref{lemma:breslau}.

Next, we may study non-local operators, where the order of differentiability is not fixed. Here is a possible example:

\begin{equation*}
\dt u - \int_\X \frac{u(y)-u(x)}{ |y-x|^{N+ \frac12 - \frac{1}{4}\sin\frac{1}{|x-y|}}}\,dy = 0\,.
\end{equation*}
We may even allow the order of differentiability to depend on $u(x), u(y)$ as in the following example:
\begin{equation*}
\dt u - \int_{B_1(x)} \frac{u(y)-u(x)}{|x-y|^{N+\frac12-\frac{1}{4}\exp(-|u(y)-u(x)|)}}\,dy = 0 \,.
\end{equation*}
The general case including precise assumptions is explained in Proposition~\ref{lemma:bielefeld}.

\subsubsection*{Potential extensions}   
Our framework could be extended and adapted to cover the following models which are not currently in its scope.

In the series of papers \cite{MR2914243, MR3218830, MR2795714, MR3286677}
the properties of solutions to the following conservation laws 
\begin{equation*}
 \dt u + \div f(u) + \mathcal{L}u = g(x,t)
\end{equation*}
have been studied. The non-local operator $\mathcal{L}$ is given by formula \eqref{L:lin}.

One may also consider  the following general fractional porous medium equation with variable density 
\begin{equation*}
\rho(x)\dt u  + (-\Delta)^s (u^{m-1}u) = 0,
\end{equation*}
which was considered in 
\cite{MR3412412, MR3158444}.

We also mention the work \cite{MR3466219}, devoted to the 
inhomogeneous non-local diffusion equation
\begin{equation*}
  \dt u (x, t) = \int_\R J\bigg(\frac{x-y}{g(y)}\bigg) \frac{u(y, t)}{g(y)}\,dy - u(x, t),
\end{equation*}
where $J$ is a non-negative even function supported
in the unit interval $[-1, 1]$ and such that  $\int_\R J(x)\,dx = 1$ and the function $g$ is continuous and positive.

 The non-linear porous medium equation with fractional potential pressure has the following form
\begin{equation}\label{pmep}
  \dt u = \div\big(u^{m_1}\nabla(-\Delta)^{-s} u^{m_2}\big)
\end{equation}
where $m_1,m_2\geq 1$ and $s\in (0,1)$.
This equation was first studied in \cite{MR2575479} in the one dimensional case and for $m_1=m_2=1$, and solutions to the corresponding Cauchy problem were shown to exist and to be unique. Moreover, an explicit self-similar compactly supported solution has been constructed for this equation for $N=1$ in~\cite{MR2575479} and for $N\geq 1$ in~\cite{MR2817383, MR3294409}. Independently, a theory of equation~\eqref{pmep}
has been developed in \cite{MR2847534,MR2773189}. Recent results and several other references have been obtained and gathered in~\cite{MR3151879, MR3419724}.

We conclude this overview by recalling the recent paper~\cite{0951-7715-32-1-1}\footnote{We learned about this model thanks to a presentation by Giuseppe Maria Coclite at a workshop at the Norwegian University of Science and Technology.}, in which the authors prove existence, uniqueness and stability of weak solutions to a non-linear, non-local, vector-valued wave equation
\begin{equation}\label{eq:wave}
 \left\{\begin{aligned}
        & \partial^2_t u(x, t) = \big(Ku(\,\cdot\,, t)\big)(x), &&x \in \R^N,\ t > 0,\\
        & u(x, 0) = u_0(x),\quad \partial_t u(x, 0) = v_0(x), &&x \in \R^N ,
        \end{aligned}
 \right.
\end{equation}
with the operator
\begin{equation}\label{eq:peridynamics}
 (Ku)(x) = \int_{B_\delta(x)} f\big(y - x, u(y) - u(x)\big)\, dy,\text{ for every $x \in \R^N$}.
\end{equation}
Here $f \in C^1(\Omega; \R^N )$ is defined on the set
$\Omega := (\R^N \setminus \{0\}) \times \R^N$ and 
\begin{equation*}
f(-y, -u) = -f(y, u)\quad\text{ for every $(y, u) \in \Omega \times \R^N$}.
\end{equation*}
Some additional regularity and integrability properties of the function $f$ are imposed, derived
from the physical assumption of hyperelasticity. In particular, the singular part
of the operator on the $(x,y)$ diagonal is assumed to be the $s$-fractional $p$-Laplacian we discussed in equation~\eqref{plap}.
This model is based on the previous results
gathered in e.g.~\cite{MR2881042}, \cite{MR3061151}, \cite{MR3297136} and it serves to study \emph{peridynamics}, a non-local elasticity theory, used to explain the formation of fractures in solids, defects, dislocations etc. While
the wave equation~\eqref{eq:wave} requires a separate set of methods to study properties of solutions, the operator defined by formula~\eqref{eq:peridynamics} has a similar structure to the one portrayed in Definition~\ref{def:jump-kernel}.

\subsection*{Outline} The paper is structured as follows.
In Section~\ref{sec:operator} we discuss the properties of the operator $\Lap$.
In Section~\ref{sec:regular} we solve equation~\eqref{pns} in a regular setting.
In Section~\ref{sec:existence} we prove our main results (see Theorem~\ref{thm:main} and Corollary~\ref{cor:main}),
namely existence and uniqueness of solutions to equation~\eqref{pns} and we study their properties.
In Section~\ref{sec:examples} we give several examples of jump kernels.
\end{section}
\begin{section}{Non-linear Lévy operator}\label{sec:operator}
We begin by introducing our notation. We use the Banach space
\begin{equation*}
\Lpp{\X}= L^1(\X)\cap L^\infty(\X) 
\end{equation*}
supplemented with the usual norm $\|u\|_{[1,\infty]} = \|u\|_1+\|u\|_\infty$ as well as standard Sobolev spaces
$W^{1,\infty}(\X)$ and $W^{1,1}_{\loc}\big([0,\infty),L^1(\X)\big)$.

  We employ the space of functions of bounded variation, following~\cite{MR3726909}. Let $u\in L^1(\X)$ 
  and for $i=1,\ldots,N$ suppose there exist finite signed Radon measures $\lambda_i$ such that
  \begin{equation*}
   \int_\X u\,\partial_{x_i} \phi\,dx = - \int_\X \phi\, d\lambda_i\quad\text{for every $\phi\in C_c^\infty(\X)$}.
  \end{equation*}
  We define
  \begin{equation*}
   |Du|(\X) = \sum_{i=1}^N\sup\bigg\{\int_\X\Phi_i\,d\lambda_i: \Phi\in C_0(\X,\X), \|\Phi\|_{C_0(\X,\X)}<1\bigg\}.
  \end{equation*}
  Then we say $u\in BV(\X)$ if the value of the following norm
  \begin{equation*}
    \|u\|_{BV} = 2\|u\|_1 + |Du|(\X)
  \end{equation*}
  is finite (caution: the number 2 in front of the $L^1$-norm is added to simplify estimates below).
  We also recall an alternative characterization of the $BV$-space which is more useful for us.
  \begin{lemma}\label{lemma:bv-char}
   Let $(e_1,\ldots,e_N)$ be the canonical basis of $\X$. Suppose $u\in L^1(\X)$ and
   denote
   \begin{equation*}
    |u|_{BV} = \sum_{i=1}^N \limsup_{h\to0^+} \int_{\X} \frac{\big|u(x)-u(x+he_i)\big|}{h}\,dx.
   \end{equation*}
   Then $u\in BV(\X)$ if and only if $|u|_{BV}<\infty$. Moreover $|u|_{BV} = |Du|(\X)$.
  \begin{proof}
  See \cite[Theorem~14.53]{MR3726909}.
  \end{proof}
  \end{lemma}
  
  \begin{lemma}\label{lemma:bvnorm}
   If $u\in BV(\X)$ then
   \begin{equation*}
    \sup_{y\in\X\setminus\{0\}}\int_\X\frac{\big|u(x)-u(x-y)\big|}{1\wedge|y|}\,dx \leq \|u\|_{BV}.
   \end{equation*}
   \begin{proof}
    According to~\cite[Lemma~14.37]{MR3726909} we have
   \begin{equation*}
    \int_\X\big|u(x)-u(x-y)\big|\,dx \leq |y|\,|Du|(\X).
   \end{equation*} 
   Thus
   \begin{equation*}
    \sup_{y\in\X\setminus\{0\}}\int_\X\frac{\big|u(x)-u(x-y)\big|}{1\wedge|y|}\,dx \leq 2\|u\|_1 + |Du|(\X) = \|u\|_{BV}.\qedhere
   \end{equation*}
   \end{proof}
  \end{lemma}

  Now we are ready to show that if $v$ is a bounded function, the linear operator $\Lap_v$ given by formula~\eqref{operator} is well-defined on the space of $BV$-functions. 
  \begin{lemma}\label{lemma:bv}
   Let $R>0$ and suppose $\meas:(\R\times\R)\times(\X\times\X)\to\R$ is a function satisfying conditions
   \ref{a1:positive} and \ref{a5:1levy} in Definition~\ref{def:jump-kernel} with a constant $K_R$.
   For every $u\in BV(\X)$ and every $v\in L^\infty(\X)$ such that $\|v\|_\infty \leq R$ we have
   \begin{equation*}
    \big[u(x)-u(y)\big]\meas_{v,x,y}\in L^1(\XX).
   \end{equation*}
   Moreover,
   \begin{equation*}
   \Lap_v:BV(\X)\to L^1(\X)\quad\text{and}\quad\|\Lap_vu\|_1 \leq K_R\|u\|_{BV}.
   \end{equation*}
   \begin{proof}
   We use conditions~\ref{a1:positive} and~\ref{a5:1levy}, as well as Lemma~\ref{lemma:bvnorm} and the Fubini-Tonelli theorem to estimate
   \begin{align*}
    \iint_{\XX}\big|u(x)-u(y)\big|\,d\meas_{v}
    &\leq\iint_{\XX}\frac{\big|u(x)-u(x-y)\big|}{1\wedge|y|}\big(1\wedge|y|\big)m_R\big(|y|\big)\,dx\,dy\\
    &\leq\|u\|_{BV}\int_{\X}\big(1\wedge|y|\big)m_R\big(|y|\big)\,dy = \|u\|_{BV}K_R.
   \end{align*}
   Then it follows that
   \begin{equation*}
    \|\Lap_vu\|_1 = \int_\X\bigg|\int_\X \big[u(x)-u(y)\big]\meas_{v,x,y}\,dy\,\bigg|\,dx \leq \|u\|_{BV}K_R.\qedhere
   \end{equation*}
   \end{proof}
  \end{lemma}
  
  Let  us now prove other properties of the operator $\Lap_v$.
  \begin{lemma}\label{lemma:sym}
    Let $R>0$ and suppose $\meas:(\R\times\R)\times(\X\times\X)\to\R$ is a function satisfying conditions \ref{a1:positive}, \ref{a2:symmetry} and
   \ref{a5:1levy} in Definition~\ref{def:jump-kernel} with a constant $K_R$.
  For every $v\in L^\infty(\X)$ such that $\|v\|_\infty \leq R$ the operator $\Lap_v:C_c^\infty(\X)\to \Lpp{\X}$ is 
   $L^2$-symmetric and positive-definite. If $u\in BV(\X)$ and $\phi\in C_c^\infty(\X)$ then
   $\int_\X\Lap_vu\phi\,dx = \int_\X\Lap_v\phi u\,dx$.
   \begin{proof}
    Let $\phi,\psi\in C_c^\infty(\X)$. We have
    \begin{equation*}
     \sup_{x\in\X}\bigg|\int_\X \big[\phi(x)-\phi(y)\big] \meas_{v,x,y}\,dy\,\bigg|\leq 2\|\phi\|_{W^{1,\infty}} K_R,
    \end{equation*}
    which combined with Lemma~\ref{lemma:bv} gives us $\Lap_v\phi\in \Lpp{\X}\subset L^2(\X)$.
    Thanks to condition~\ref{a2:symmetry} in Definition~\ref{def:jump-kernel} we may ``symmetrize'' 
    the double integral and obtain (see Remark~\ref{rem:sym} below)
    \begin{multline*}
     \int_\X\psi(x)\big(\Lap_v\phi\big)(x)\,dx 
     = \int_\X\psi(x)\int_\X \big[\phi(x)-\phi(y)\big] \meas_{v,x,y}\,dy\,dx\\
     = \frac{1}{2}\iint_{\XX} \big[\psi(x)-\psi(y)\big]\big[\phi(x)-\phi(y)\big] \,d\meas_{v}
     = \int_\X\big(\Lap_v\psi(x)\big)\phi(x)\,dx.
    \end{multline*}
    The same observation holds if we exchange $\psi$ by $u\in BV(\X)\subset L^1(\X)$ (all
    the integrals are well-defined because of Lemma~\ref{lemma:bv}).
    Finally
    \begin{equation*}
     \int_\X\phi(x)\big(\Lap_v \phi\big)(x)\,dx =
     \iint_{\XX}\big[\phi(x)-\phi(y)\big]^2 \,d\meas_{v}\geq 0.\qedhere
    \end{equation*}
   \end{proof}
  \end{lemma}
  \begin{remark}\label{rem:sym}
   The ``symmetrization argument'' that we use in the proof of Lemma~\ref{lemma:sym} is going the be used
   extensively throughout this paper. In every instance it looks like the following identity
   \begin{multline*}
     \iint_{\XX} \big[f(x)-f(y)\big] g(x)\meas_{v,x,y}\,dy\,dx
     = \iint_{\XX} \big[f(y)-f(x)\big] g(y)\meas_{v,y,x}\,dx\,dy\\
     = \frac{1}{2}\iint_{\XX} \big[f(x)-f(y)\big]\big[g(x)-g(y)\big] d\meas_{v,x,y}\,dy\,dx,
   \end{multline*}
   for functions $f$, $v$ and $g$ such that 
   \begin{equation*}
    \big[f(x)-f(y)\big]\meas_{v,x,y}\in L^1(\XX)\quad\text{and}\quad g\in L^\infty(\X)
   \end{equation*}
   (notably we may take $g\equiv 1$ and the integral vanishes). The first equality requires no effort, it is simply
   renaming the variables. In the second equality we use the Fubini-Tonelli theorem to exchange $dx$ and $dy$, 
   use the symmetry property~\ref{a2:symmetry} of the jump kernel introduced in Definition~\ref{def:jump-kernel},
   which states that $\meas_{v,y,x}=\meas_{v,x,y}$, and take the average of both integrals.
  \end{remark}

  In the next theorem, we prove a result which is reminiscent of the famous Kato inequality~\cite{MR0333833,MR2056467}.
  In a way, it is central to our entire work and a source of some of its limitations.
  In particular, it is the reason behind the monotonicity condition~\ref{a3:monotonicity}
  in Definition~\ref{def:jump-kernel} as well as the introduction of the $BV$ space.
  Indeed, in Lemma~\ref{lemma:bvestimate} we show that the solution $u(t)$ belongs to $BV(\X)$ for all $t> 0$ if the initial
  condition does too. Then we we may use the theorem thanks to Lemma~\ref{lemma:bv}.
  \begin{theorem}\label{thm:kato}
   Let $\meas:(\R\times\R)\times(\X\times\X)\to\R$ be a function satisfying conditions \ref{a1:positive}, \ref{a2:symmetry} and
   \ref{a3:monotonicity} in Definition~\ref{def:jump-kernel}.
   If $u,v\in L^\infty(\X)$ are such that
   \begin{equation*}
   \big[u(x)-u(y)\big]\meas_{u,x,y}\in L^1(\XX),\quad \big[v(x)-v(y)\big]\meas_{v,x,y}\in L^1(\XX)
   \end{equation*}
   then
   \begin{equation*}
    \int_\X(\Lap_uu-\Lap_vv)\sgn (u-v)\,dx\geq 0.
   \end{equation*}
  \begin{proof}
      Let 
      \begin{align*}
       &\eta(x) = \sgn\big(u(x)-v(x)\big),\\
       &f(x,y) = \big[u(x)-u(y)\big]\meas_{u,x,y}-\big[v(x)-v(y)\big]\meas_{v,x,y}.
      \end{align*}
      Since we assume $\Lap_uu-\Lap_vv\in L^1(\X)$ and since $\eta\in 
L^\infty(\X)$, the symmetrization argument (see Remark~\ref{rem:sym}) 
gives us
    \begin{multline*}
      \int_\X \big((\Lap_uu)(x)-(\Lap_vv)(x)\big) \sg{x}\,dx\\
       =\iint_{\XX} f(x,y) \,dy\,\sg{x}\,dx 
       =\frac{1}{2}\iint_{\XX} f(x,y)\big[\sg{x}-\sg{y}\big]\,dy\,dx.
    \end{multline*}
Because of the symmetry condition~\ref{a2:symmetry}, we have $f(x,y) = 
-f(y,x)$ and hence
\begin{align*}
&\frac12 \iint_{\XX} f(x,y)\big[\sg{x}-\sg{y}\big]\,dy\,dx \\
&\qquad=\frac12 \iint_{\XX}f(x,y)\big[\sg{x}-\sg{y}\big] \mathbbm{1}_{\{\sg{x} > \sg{y}\}} \,dy\,dx \\
&\qquad\quad+ \frac12 \iint_{\XX} f(x,y)\big[\sg{x}-\sg{y}\big] \mathbbm{1}_{\{\sg{x} < \sg{y}\}} \,dy\,dx
\\
&= \iint_{\XX} f(x,y)\big[\sg{x}-\sg{y}\big] \mathbbm{1}_{\{\sg{x} > \sg{y}\}} 
\,dy\,dx \,.
\end{align*}
The set $\{\sg{x} > 
\sg{y}\}$ is a subset of the set $M=\{ u(x) \geq v(x), \; v(y) \geq u(y) \}$. 

In order to complete the proof, if suffices to show that $f$ is non-negative on $M$. If 
$(x,y) \in M$ and additionally $v(x) \geq v(y)$, then $f(x,y) \geq 0$ because of 
assumption~\ref{a3:monotonicity}. Thus, we only need to consider the 
situation, where $(x,y) \in M$ and $v(x) < v(y)$. There are two cases. If 
$u(x) \geq u(y)$, then $f(x,y) \geq 0$ just because $\meas$ is non-negative. If 
$u(x) < u(y)$, then $ v(x) \leq u(x) < u(y) \leq v(y)$. In this case,  
assumption~\ref{a3:monotonicity} implies
\begin{align*}
\big[v(y)-v(x)\big] \meas(v(y), v(x), 
y,x)-\big[u(y)-u(x)\big] \meas(u(y), u(x),y,x) \geq 0.
\end{align*}
Due to assumption~\ref{a2:symmetry} we 
obtain 
\begin{align*}
\big[u(x)-u(y)\big] \meas(u(x), u(y), 
x,y)-\big[v(x)-v(y)\big] \meas(v(x), v(y),x,y) \geq 0\,,
\end{align*}
which is nothing but $f(x,y) \geq 0$. The proof of Theorem \ref{thm:kato} is 
complete. 
\end{proof}
  \end{theorem}
\end{section}
\begin{section}{Regular jump kernels}\label{sec:regular}

In this section we construct global-in-time unique classical solutions of 
problem~\eqref{pns} under strong regularity assumptions on the jump kernel. 
Notice that we do not assume the jump kernel to be homogeneous.

  \begin{definition}\label{def:regular}
  We say a function $\meas:(\R\times\R)\times(\X\times\X)\to\R$ is a regular jump kernel
  if it satisfies conditions
  \ref{a1:positive}, \ref{a2:symmetry} and \ref{a3:monotonicity} in Definition~\ref{def:jump-kernel} and in addition:
  \begin{description}
    \item[(B1)\label{b1:nonsing}] it is integrable on the diagonal $y=x$, namely, 
          for each $R>0$ there exists $M_R>0$ such that 
          \begin{equation*}\hspace*{\leftmargin}
            \sup_{\|u\|_\infty\leq R}\sup_{x\in\X} \int_\X \meas_{u,x,y}\,dy = M_R<\infty;
          \end{equation*}
    \item[(B2)\label{b2:Lip:reg}] it is locally Lipschitz-continuous with respect to $u$, that is, for each $R>0$ 
    there exists $L_R>0$ such that 
          \begin{equation*}\hspace*{\leftmargin}
            \sup_{x\in\X} \int_\X |\meas_{u,x,y}-\meas_{v,x,y}|\,dy
            \leq L_R \|u-v\|_{[1,\infty]}
          \end{equation*}
          for all $\|u\|_\infty\leq R$ and $\|v\|_\infty\leq R$.
  \end{description}
  \end{definition}
  
  First we prove a counterpart of Lemma~\ref{lemma:bv} for regular jump kernels.
  It turns out that in this case for every $v\in\Lpp{\X}$ the linear operator $\Lap_v$ is bounded on $\Lpp{\X}$.
\begin{lemma}\label{lemma:bounded}
 If $\meas$ is a regular jump kernel then for every $u,v\in\Lpp{\X}$ we have $\Lap_vu\in\Lpp{\X}$.
  \begin{proof}
    Let $R$ be such that $\|v\|_\infty\leq R$. It follows from condition~\ref{b1:nonsing}
    that 
    \begin{equation*}
    \|\Lap_vu\|_\infty \leq 2 \|u\|_\infty M_R.
    \end{equation*}
    Renaming the variables, using the symmetry condition~\ref{a2:symmetry} (cf.~Remark~\ref{rem:sym}) and property~\ref{b1:nonsing} we obtain
    \begin{equation*}
      \iint_{\XX} \big|u(y)\big| \meas_{v,x,y}\,dy\,dx = \iint_{\XX} \big|u(x)\big| \meas_{v,x,y}\,dy\,dx
      \leq \|u\|_1 M_R
    \end{equation*}
    and therefore 
    \begin{equation*} 
      \|\Lap_vu\|_1 \leq \iint_{\XX} \big(\big|u(x)\big|+\big|u(y)\big|\big)\meas_{v,x,y}\,dy\,dx\leq 2 \|u\|_1 M_R,
    \end{equation*}
    which completes the proof that  $\Lap_v:\Lpp{\X}\to\Lpp{\X}$.
 \end{proof}
\end{lemma}
  \begin{lemma}\label{lemma:regular-lipschitz}
   If $\meas$ is a regular jump kernel then the operator $F(u)=-\Lap_uu$
    is locally Lipschitz as a mapping $F: \Lpp{\X} \to \Lpp{\X}$.
   \begin{proof}
    Let $u,v\in \Lpp{\X}$ be such that $\|u\|_\infty\leq R$ and $\|v\|_\infty\leq R$. 
    Note the identity
    \begin{equation*}
      \Lap_uu-\Lap_vv= \Lap_u(u-v)+(\Lap_u-\Lap_v)v.
    \end{equation*}
    We have, using the integrability condition~\ref{b1:nonsing},
    \begin{multline*} 
        \|\Lap_u(u-v)\|_\infty 
        = \sup_{x\in\X}\bigg|\int_\X \big(u(x)-u(y)-v(x)+v(y)\big)\meas_{u,x,y}\,dy\,\bigg|\\
        \leq 2\|u-v\|_\infty\sup_{x\in\X}\int_\X\meas_{u,x,y}\,dy \leq 2M_R\|u-v\|_\infty
    \end{multline*}
    and by the local Lipschitz-continuity of the jump kernel~\ref{b2:Lip:reg} we obtain
    \begin{multline*} 
        \|\Lap_uv-\Lap_vv\|_\infty =
        \sup_{x\in\X}\bigg|\int_\X\big[v(x)-v(y)\big]\big(\meas_{u,x,y}-\meas_{v,x,y}\big)\,dy\,\bigg|\\
        \leq 2\|v\|_\infty\sup_{x\in\X}\int_\X|\meas_{u,x,y}-\meas_{v,x,y}|\,dy 
        \leq 2L_R\|v\|_\infty\|u-v\|_{[1,\infty]}.
    \end{multline*}
    By a similar calculation and the symmetrization argument (see Remark~\ref{rem:sym}) we also get
    \begin{align*}
      \|\Lap_u(u-v)\|_1
      &= \int_\X\bigg|\int_\X \big[u(x)-u(y)-v(x)+v(y)\big]\meas_{u,x,y}\,dy\,\bigg|\,dx\\
      &\leq \iint_{\XX} \big(\big|u(x)-v(x)\big|+\big|v(y)-u(y)\big|\big)\meas_{u,x,y}\,dy\,dx\\
      &= 2\int_\X \big|u(x)-v(x)\big|\int_\X\meas_{u,x,y}\,dy\,dx
      \leq 2M_R\|u-v\|_1\
    \end{align*}
    and
    \begin{align*}
      \|\Lap_uv&-\Lap_vv\|_1
      = \int_\X\bigg|\int_\X\big[v(x)-v(y)\big](\meas_{u,x,y}-\meas_{v,x,y})\,dy\,\bigg|\,dx\\
      &\leq \iint_{\XX} \big(\big|v(x)\big|+\big|v(y)\big|\big)|\meas_{u,x,y}-\meas_{v,x,y}|\,dy\,dx\\
      &= 2\int_\X\big|v(x)\big|\int_\X|\meas_{u,x,y}-\meas_{v,x,y}|\,dy\,dx
      \leq 2L_R\|v\|_1\|u-v\|_{[1,\infty]},
    \end{align*}
    which completes the proof of Lemma~\ref{lemma:regular-lipschitz}.
   \end{proof}
  \end{lemma}

  Now we may construct local-in-time solutions via the Banach fixed point argument.  
  
  \begin{lemma}\label{lemma:loc-clas}
   If $\meas$ is a regular jump kernel then for every $u_0\in\Lpp\X$ there exist $T>0$ and a unique 
   local classical solution $u\in C^1\big([0,T],\Lpp{\X}\big)$ to problem~\eqref{pns} on $[0,T]$.
   \begin{proof}
    Notice that if $v\in C^1\big([0,T],\Lpp{\X}\big)$ then the expression $\dt v +\Lap_vv$ is well-defined
    for every regular jump kernel $\meas$ (see Lemma~\ref{lemma:bounded}).
    Consider the mapping $F(v)=-\Lap_vv$ and an integral operator
    \begin{equation*}
      \mathfrak{F}v(t)=u_0+\int_0^t F\big(v(s)\big)\;ds
    \end{equation*}
    in the Banach space $C\big([0,T], \Lpp{\X}\big)$.
    We know from Lemma~\ref{lemma:regular-lipschitz} that the operator $F$ is locally Lipschitz.
    Therefore it suffices to apply the Banach contraction principle on a certain interval $[0,T]$ in order
    to obtain the unique fixed point $u$ of the operator $\mathfrak{F}$. Moreover, 
    \begin{equation*}
      \mathfrak{F}:C\big([0,T], \Lpp{\X}\big)\to C^1\big([0,T], \Lpp{\X}\big),               
    \end{equation*}
   hence  the equation $\dt u = F(u)= -\Lap_uu$ is satisfied in the classical sense.
   \end{proof}
  \end{lemma}
  \begin{lemma}\label{lemma:Lp}
    If $u$ is a local classical solution to problem~\eqref{pns}	on $[0,T]$ 
    then for every $t\in[0,T]$ and every $p\in[1,\infty]$
    we have
    \begin{equation}\label{Lp:ns}
      \|u(t)\|_p\leq \|u_0\|_p.
    \end{equation}
  \begin{proof}
    We know that $u\in C^1\big([0,T],L^p(\X)\big)$ for every $p\in[1,\infty]$.
    Let us fix $p\in(1,\infty)$.
    We multiply the equation 
    $\dt u = -\Lap_uu$
    by $|u|^{p-2}u$ and integrate  with respect to $x$ to obtain 
    \begin{equation*}
      \int_\X \big(\dt u(x)\big)\big(\big|u(x)\big|^{p-2}u(x)\big)\,dx 
      = \iint_{\XX}\big[u(y)-u(x)\big]\big|u(x)\big|^{p-2}u(x)\,d\meas_u.
    \end{equation*}
    We obtain
    \begin{equation*}
      \dt\int_\X \big|u(x)\big|^p\,dx = p\int_\X \big(\dt u(x)\big)\big(\big|u(x)\big|^{p-2}u(x)\big)\,dx
    \end{equation*}
    and thanks to the symmetrization argument (see Remark~\ref{rem:sym}) it follows that
    \begin{multline*}
      \iint_{\XX}\big[u(y)-u(x)\big]\big|u(x)\big|^{p-2}u(x)\,d\meas_u
      \\= \frac{1}{2}\iint_{\XX} \big[u(y)-u(x)\big]\big(\big|u(x)\big|^{p-2}u(x)
      -\big|u(y)\big|^{p-2}u(y)\big)\,d\meas_u\leq 0.
    \end{multline*}
    The last inequality holds because the mapping $a\mapsto |a|^{p-2}a$
    is non-decreasing on $\mathbb{R}$ and the measure $d\meas_u$ is non-negative due to condition~\ref{a1:positive}.
    We thus have proved
    inequalities~\eqref{Lp:ns} for all $p\in (1,\infty)$.
    Because the function $p\mapsto \|f\|_p$ is continuous for every $f\in\Lpp{\X}$,
    we may pass to the limits with $p\to 1$ and $p\to\infty$ to obtain inequalities in~\eqref{Lp:ns} for all $p\in[1,\infty]$.
  \end{proof}
  \end{lemma}
  
  Now we are ready to prove existence of solutions in the case of regular jump kernels.
  
  \begin{theorem}\label{thm:classical}
  If $\meas$ is a regular jump kernel then the classical solution is global.
  \begin{proof}
    Consider the local classical solution $\dt u = -\Lap_uu$ on an interval 
    $[0,T]$, as constructed in Lemma~\ref{lemma:loc-clas}.
    It follows from Lemma~\ref{lemma:Lp} that
    \begin{equation*}
    \|u(t)\|_{[1,\infty]}\leq\|u_0\|_{[1,\infty]}
    \end{equation*}
    thus the local classical solution
    may be extended to all $t\in [0,\infty)$ by a usual continuation 
    argument.
    \end{proof}
  \end{theorem}
    We now examine some of the properties of classical solutions, which will be useful in the next section.
    Notice that these results cannot be directly applied in the general case, where we need a weaker
    notion of solutions.
    In the following lemma we discuss the $L^1$-contraction property.
    In the proof we use the Kato-type inequality from Theorem~\ref{thm:kato}.
    \begin{lemma}[$L^1$-contraction for classical solutions]\label{lemma:contraction}
    If $u,v$ are classical solutions to problem \eqref{pns} with initial conditions $u_0$ and $v_0$, respectively
    then
    \begin{equation*}
       \|u(t)-v(t)\|_1 \leq \|u_0-v_0\|_1
      \end{equation*}
      for every $t\geq 0$.
    \begin{proof}
    We have $u,v\in C^1\big([0,\infty),\Lpp\X\big)$, therefore the following integral
    \begin{equation*}
     \int_0^\infty\int_\X \Big(\dt \big(u(t,x)-v(t,x)\big) +\big(\Lap_uu-\Lap_vv\big)(t,x)\Big)\psi(t,x)\,dx\,dt = 0
    \end{equation*}
    is convergent for 
    \begin{equation*}
     \psi(t,x) = \mathbbm{1}_{[0,T]}(t)\sgn\big((u(t,x)-v(t,x)\big),
    \end{equation*}
    where we arbitrarily fix $T>0$.
    Thus by the assumed regularity of $u$ and $v$ we get
    \begin{equation*}
     \int_0^T\dt  \int_\X |u-v|\,dx\,dt
     = -\int_0^T\int_\X\Big(\Lap_uu-\Lap_vv\Big)\sgn(u-v)\,dx\,dt.
    \end{equation*}
    It follows from Theorem~\ref{thm:kato} that
    \begin{equation*}
    \int_0^T\dt \|u-v\|_1\,dt \leq 0
    \end{equation*}
    and consequently $\|u(T)-v(T)\|_1 \leq \|u_0-v_0\|_1$.
  \end{proof}
  \end{lemma}
    In the next lemma we estimate the $BV$-norm of a solution in case of an additional
    assumption of homogeneity (as in condition~\ref{a4:homogeneity} in Definition~\ref{def:jump-kernel}).
    This estimate will help us to establish relative compactness of an approximating sequence of solutions
    we construct in Lemma~\ref{lemma:approx} by regularizing the jump kernel. 
  \begin{lemma}\label{lemma:bvestimate}
   Let $\meas$ be a regular jump kernel satisfying the homogeneity condition~\ref{a4:homogeneity} i.e.
   \begin{equation*}\meas\big(v(x),v(y);x,y\big) = \meas\big(v(x),v(y);|x-y|\big).\end{equation*}
   If $u_0\in BV(\X)$ and $u\in C^1\big([0,T],\Lpp{\X}\big)$ is the classical solution to problem~\eqref{pns}
   then
   \begin{equation*}
    \|u(t)\|_{BV} \leq \|u_0\|_{BV}\quad\text{for every $t\geq0$}.
   \end{equation*}
   \begin{proof}
    Let $v_\xi(x)=v(x+\xi)$ for an arbitrary $\xi\in\X$ and $v\in L^\infty(\X)$. 
    Because $\meas$ is homogeneous
    we have $\meas_{v,x+\xi,y+\xi}=\meas_{v_\xi,x,y}$. Consequently, for an arbitrary $w\in L^\infty(\X)$ we have
    \begin{multline}\label{eq:shift}
     \big(\Lap_v w\big)(x+\xi)=\int_\X \big(w(x+\xi)-w(y)\big)\meas_{v,x+\xi,y}\,dy \\=
     \int_\X \big(w_\xi(x)-w_\xi(y)\big)\meas_{v_\xi,x,y}\,dy =
     \big(\Lap_{v_\xi}w_\xi\big)(x),
    \end{multline}
    and the integrals are convergent because $\meas$ satisfies condition~\ref{b1:nonsing}.
    
    It follows from identity~\eqref{eq:shift}, and assumed regularity of $u$, that $u_\xi$ is the classical solution to problem~\eqref{pns} with initial condition $u_{0,\xi}(x) = u_0(x+\xi)$.
    By Lemma~\ref{lemma:contraction} we thus have
    \begin{equation*}
     \|u_\xi(t)-u(t)\|_1 \leq \|u_{0,\xi}-u_0\|_1
    \end{equation*}
      and by \cite[Lemma~14.37]{MR3726909}
    \begin{equation*}
     \|u_{\xi,0}-u_0\|_1 = \int_\X\big|u_0(x+\xi)-u_0(x)\big|\,dx \leq |\xi|\|u_0\|_{BV}.
    \end{equation*}
    By taking $\xi=he_i$ with $\{e_i\}$ being the canonical basis of $\X$ and $h>0$ we get
    \begin{equation*}
     \int_\X\frac{\big|u(x)-u(x+he_i)\big|}{h}\,dx \leq \|u_0\|_{BV}
    \end{equation*}
    and it follows by Lemma~\ref{lemma:bv-char} that $\|u(t)\|_{BV} \leq \|u_0\|_{BV}$.
   \end{proof}
  \end{lemma}
  \end{section}
  \begin{section}{Construction of solutions}\label{sec:existence}
\subsection*{Definitions of solutions}
Our goal is to construct solutions in the case of general homogeneous jump kernels and  prove the results stated in Theorem~\ref{thm:main}.
As usual in such contexts we work with a weak formulation of problem~\eqref{pns}.
\begin{definition}\label{def:distributional}
    We say $u\in L^\infty\big([0,\infty)\times\X\big)$ is a very weak solution to problem~\eqref{pns} if 
    \begin{equation}\label{pns:weak}
      \int_0^\infty\int_\X u(t,x)\big[\big(\dt-\Lap_u\big)\psi\big](t,x)\,dx\,dt = 0
    \end{equation}
    for every $\psi\in C^\infty_c\big((0,\infty)\times \X\big)$ and 
    $\lim_{t\to0}u(t,x) = u_0(x)$ in $L^1_{\loc}(\X)$.
\end{definition}
\begin{remark}
 Notice that thanks to Lemma~\ref{lemma:bv} this definition is well-posed.
\end{remark}
\begin{definition}\label{def:bvsolution}
  We say a very weak solution to problem~\eqref{pns} $u$ is a strong solution to problem~\eqref{pns} if
   \begin{equation*}
    u\in L^\infty\big([0,\infty),BV(\X)\cap L^\infty(\X)\big)\cap C\big([0,\infty),L^1_{\loc}(\X)\big).
   \end{equation*}
\end{definition}
\begin{remark}\label{rem:classical-bv}
 Notice that because of Lemmas~\ref{lemma:bv} and~\ref{lemma:sym}, if $u$ is a strong solution to problem~\eqref{pns} then
 $\Lap_uu$ is well-defined and we have 
    \begin{equation*}
      \int_0^\infty\int_\X u(t,x)\dt\psi(t,x)-\big(\Lap_uu\big)(t,x)\psi(t,x)\,dx\,dt = 0.
    \end{equation*}
\end{remark}
  \begin{remark}\label{rem:constant}
   Notice that each constant function $u(t,x)\equiv k\in\R$ is
   a classical solution to problem~\eqref{pns} for every jump kernel $\meas$
   and we have $\Lap_uu \equiv 0$.
  \end{remark}
\subsection*{Approximate solutions}
In this part we regularize an arbitrary homogeneous jump kernel and study compactness of the corresponding sequence of approximations. 
In order to simplify our reasoning we begin with the following observation.
\begin{remark}\label{rem:lip-sym}
    Consider a homogeneous jump kernel $\meas$.
    For arbitrary $a,b,c,d\in \R$ by the symmetry assumption~\ref{a2:symmetry}
    we have 
    \begin{multline*}
    \big| \meas(a,b;x,y)-\meas(c,d;x,y)\big| \\
    \leq  \big| \meas(a,b;x,y)-\meas(c,b;x,y)\big| + \big| \meas(b,c;x,y)-\meas(d,c;x,y)\big|.
    \end{multline*}
    Therefore it is sufficient (and necessary) to verify the Lipschitz-continuity part of condition~\ref{a6:lipschitz}
    only for the difference
    \begin{equation*}
     \big| \meas(a,b;x,y)-\meas(c,b;x,y)\big|
    \end{equation*}
    and $a,b,c\in[-R, R]$ such that $|a-b|\geq\epsilon$ and $|c-b|\geq\epsilon$.
  \end{remark}

  \begin{lemma}\label{lemma:approx}
   For every $\epsilon\in(0,1]$
   consider a function $h_\epsilon\in C^\infty\big([0,\infty)\big)$ which is non-decreasing and such that $h_\epsilon(x) = 0$ for 
   $x \leq \frac{\epsilon}{2}$ and $h_\epsilon(x) =1$ for $x\geq \epsilon$.
   Let $\meas$ be a homogeneous jump kernel and
   \begin{equation*}
      \kernelel{a}{b} 
      = h_\epsilon\big(|a-b|\big)\mathbbm{1}_{|x-y|\geq \epsilon}(x,y)\kernell{a}{b}
    \end{equation*}
    Then $\mease$ are regular, homogeneous jump kernels.
    \begin{proof}
    Conditions~\ref{a1:positive} to~\ref{a6:lipschitz} in Definition~\ref{def:jump-kernel} are easy to verify.
    Let us check conditions~\ref{b1:nonsing} and~\ref{b2:Lip:reg} of Definition~\ref{def:regular}. For an arbitrary $R>0$ consider a pair of functions $u,v\in L^\infty(\X)$ such that $\|u\|_\infty,\|v\|_\infty\leq R$.
    Because $\meas$ is a homogeneous jump kernel, we have
    \begin{multline*}
      \sup_{x\in\X}\int_\X \kernelel{a}{b}\,dy \leq \int_\X \mathbbm{1}_{|y|>\epsilon}(y)m_R\big(|y|\big)\,dy\\
      \leq  \epsilon^{-1}\int_\X \big(1\wedge|y|\big) m_{R}\big(|y|\big)\,dy
      = \epsilon^{-1}K_{R},
    \end{multline*}
    where $m_R$ and $K_R$ come from condition~\ref{a5:1levy} satisfied by the jump kernel $\meas$.
    This verifies condition~\ref{b1:nonsing}.
    
    Notice that because of the properties of the function $h_\epsilon$ and condition~\ref{a6:lipschitz}, the jump kernel $\meas_\epsilon$ is
    locally Lipschitz-continuous in the first two variables, including the diagonal. Namely,
    \begin{equation*}
     \big|\kernelel{a}{b}-\kernelel{c}{b}\big|\leq C_{R}|a-c|\mathbbm{1}_{|x-y|>\epsilon}(x,y)m_R\big(|x-y|\big)
    \end{equation*}
    for all $-R\leq a,b,c\leq R$, where $C_R=C_{R,\frac{\epsilon}{2}}$ is the constant from condition~\ref{a6:lipschitz} satisfied by 
    the jump kernel $\meas$ (cf. Remark~\ref{rem:lip-sym}). Therefore
    \begin{align*}
      &\sup_{x\in\X} \int_\X |\mease_{u,x,y}-\mease_{v,x,y}|\,dy\\
      &\leq\sup_{x\in\X} \int_\X \big|\mease_{u,x,y}-\mease\big(u(x),v(y),x,y\big)\big|
      + \big|\mease\big(u(x),v(y),x,y\big)-\mease_{v,x,y}\big|\,dy\\
      &\leq C_R \sup_{x\in\X} \int_\X \Big(\big|u(y)-v(y)\big| + \big|u(x)-v(x)\big|\Big)\mathbbm{1}_{|x-y|>\epsilon}m_R\big(|x-y|\big)\,dy\\
      &\leq L_R \|u-v\|_{[1,\infty]},
    \end{align*}
    which confirms condition~\ref{b2:Lip:reg}.
    \end{proof}
  \end{lemma}

  \begin{theorem}[Existence of strong solutions]\label{thm:existence}
    If $\meas$ is a homogeneous jump kernel then there exists a strong solution
    to problem~\eqref{pns} for every initial condition $u_0\in BV(\X)\cap L^\infty(\X)$.
  \begin{proof}
    For every $\epsilon\in (0,1]$ consider the unique classical solution 
    $u^\epsilon$ of the following initial value problem
    \begin{equation}\label{eq:approx}
    \left\{
       \begin{aligned}
        &\partial_t  u^\epsilon  + \Lap^\epsilon_{u^\epsilon}u^\epsilon = 0,\\
        & \Lap^\epsilon_vu = \int_\X \big[u(x)-u(y)\big]\mease\big(v(x),v(y);x,y\big)\,dy,\\
        &u^\epsilon(0,\,\cdot\,)=u_0,
       \end{aligned}
     \right.
    \end{equation}
    where the regular, homogeneous jump kernels $\mease$ are introduced in Lemma~\ref{lemma:approx}.
    
    Let $R$ be a constant such that $\|u_0\|_\infty \leq R$.
    Thanks to Lemma~\ref{lemma:Lp} we have $\|u^\epsilon(t)\|_\infty \leq \|u_0\|_\infty \leq R$.
    It follows from Lemmas~\ref{lemma:bvestimate} and~\ref{lemma:bv} that 
    $\Lap^\epsilon_{u^\epsilon(t)}u^\epsilon(t)\in L^1(\X)$ is
    well-defined for every $t\geq 0$. Moreover, we get an estimate on the time derivative
    \begin{equation*}
    \big\|\dt u^\epsilon(t)\big\|_1 = \big\|\Lap^{\epsilon}_{u^{\epsilon}(t)}u^{\epsilon}(t)\big\|_1 \leq \big\|u^{\epsilon}(t)\big\|_{BV}K_R \leq \|u_0\|_{BV}K_R.
    \end{equation*}    
    For an arbitrary bounded open set $\Omega\subset\X$ with a sufficiently regular boundary
    we have the compact embedding $BV(\Omega)\subset L^1(\Omega)$ 
    (see~\cite[Theorem~14.39]{MR3726909}).
    Thus we obtain a convergent subsequence in the usual way,
    by applying the Aubin-Lions-Simon lemma (see~\cite[Theorem 1]{MR916688})
    in the space $L^\infty\big([0,T],L^1(\Omega)\big)$ for an arbitrarily fixed $T>0$. Namely, it follows that there exists a function $u$ and a subsequence $\{\epsilon_j\}$
    such that 
    \begin{equation*}
     u^{\epsilon_j}\to u\quad\text{in}\quad C\big([0,T],L^1_{\loc}(\X)\big).
    \end{equation*}
 In particular we also have $\lim_{j\to\infty}u^{\epsilon_j}(t,x)= u(t,x)$ almost everywhere. We may enhance this result by applying the Fatou lemma to inequalities
    \begin{equation}\label{p:est}
     \|u^{\epsilon}(t)\|_p \leq \|u_0\|_p,
    \end{equation}
    which we proved in Lemma~\ref{lemma:Lp}, to obtain the same estimates for $\|u\|_p$.
    By using~\cite[Theorem~14.39]{MR3726909} and Lemma~\ref{lemma:bv-char} we
    also get
    \begin{equation}\label{BV:est}
    \|u(t)\|_{BV} \leq \|u_0\|_{BV}.
    \end{equation}
    In this way we obtain that 
    \begin{equation}\label{eq:estimates}
    u\in L^\infty\big([0,\infty),BV(\X)\cap L^\infty(\X)\big)\cap C\big([0,\infty),L^1_{\loc}(\X)\big).
    \end{equation}
    Consider the sequence of approximations $u^j = u^{\epsilon_j}$. Let
    \begin{equation*}
     f_j(t,x,y) = \big[u^{j}(t,x)-u^{j}(t,x-y)\big]\meas^{\epsilon_j}\big(u^j(x),u^j(x-y),|y|\big)\psi(t,x),
    \end{equation*}
    with an arbitrary test function $\psi\in C^\infty_c\big((0,\infty)\times \X\big)$.
    Let $f$ be defined analogously for the function $u$ in place of $u^j$ and the jump kernel $\meas$
    instead of $\meas^{\epsilon_j}$.
     We have (see Remark~\ref{rem:classical-bv})
    \begin{equation}\label{eq:distributional-classical}
    \int_0^\infty\int_\X u^j\,\dt\psi\,dx\,dt
    =\int_0^\infty\iint_\XX f_j(t,x,y)\,dx\,dy\,dt
   \end{equation}
    for every $j\in\mathbb{N}$.
    
    On the left-hand side of equalities~\eqref{eq:distributional-classical} we have 
    \begin{equation*}
     \big|u^j(t,x)\dt\psi(t,x)\big| \leq \|u_0\|_\infty\big|\dt\psi(t,x)\big|,
    \end{equation*}
    hence we may pass to the limit by the Lebesgue dominated convergence theorem.
    On the right-hand side of equalities~\eqref{eq:distributional-classical}, 
    by the integrability condition~\ref{a5:1levy} satisfied by the jump kernel $\meas$,
    we obtain the following estimate
    \begin{equation}\label{eq:majorant1}
     \big|f_j(t,x,y)\big|
     \leq 2R\, m_R\big(|y|\big)\, |\psi(t,x)|.
    \end{equation}
    Because of the continuity of the jump kernel assumed in condition~\ref{a6:lipschitz} we also have
    \begin{equation*}
     \lim_{j\to\infty} f_j(t,x,y)= f(t,x,y)\quad\text{a.e.~in $(t,x,y)$}    
    \end{equation*}
    and therefore we may pass to the limit by the Lebesgue dominated convergence theorem and get
    \begin{equation*}
     \lim_{j\to\infty}\int_\X f_j(t,x,y)\,dx = \int_\X f(t,x,y)\,dx,\quad \text{a.e.~in $(t,y)$}.
    \end{equation*}
    Then, by Lemmas~\ref{lemma:bvnorm} and~\ref{lemma:bvestimate}, we have
    \begin{equation}\label{eq:majorant2}
    \begin{split}
    \bigg|&\int_\X f_j(t,x,y)\,dx\,\bigg|\\
    &\leq\sup_{x\in\X}|\psi(t,x)|\big(1\wedge|y|\big)m_R\big(|y|\big) \sup_{y\in\X\setminus\{0\}}\int_\X \frac{\big|u^{j}(t,x)-u^{j}(t,x-y)\big|}{1\wedge|y|}\,dx
    \\
    &\leq \sup_{x\in\X}|\psi(t,x)|\big(1\wedge|y|\big)m_R\big(|y|\big)\|u_0\|_{BV}.
    \end{split}
    \end{equation}
 By the Lebesgue dominated convergence theorem we may thus pass to the limit once more. Combining both arguments
 in~\eqref{eq:majorant1} and~\eqref{eq:majorant2}, we get
 \begin{equation*}
  \lim_{j\to\infty}\iint_\XX f_j(t,x,y)\,dx\,dy
  = \iint_\XX f(t,x,y)\,dx\,dy\quad \text{a.e.~in $t$}.
 \end{equation*}
 We also have 
 \begin{equation*}
  \bigg|\iint_\XX f_j(t,x,y)\,dx\,\bigg| \leq \mathbbm{1}_{\supp\psi}(t)\|\psi\|_\infty \|u_0\|_{BV}K_R,
 \end{equation*}
 which allows us to pass to the limit with the integral in time.
 In this way we have shown that $u$ satisfies the following integral 
 equality
 \begin{equation*}
   \int_0^\infty\int_\X u\,\dt\psi\,dx\,dt
    =\int_0^\infty\iint_\XX \big[u(t,x)-u(t,y)\big]\meas_{u,x,y}\psi(t,x)\,dx\,dy\,dt
 \end{equation*}
 for each test function $\psi$.
 To finish the proof we apply the Fubini-Tonelli theorem, which is possible thanks to
 Lemma~\ref{lemma:bv}, to exchange $dx$ and $dy$.
 This together with~\eqref{eq:estimates} confirms Definition~\ref{def:distributional}.
   \end{proof}
  \end{theorem}
  \begin{theorem}\label{thm:regularity}
   If $u$ is a strong solution to problem~\eqref{pns}
   then \begin{equation*}u\in W^{1,1}_{\loc}\big([0,\infty),L^1(\X)\big).\end{equation*}
   \begin{proof}
  Consider a sequence $\{\phi_n\}\in C_c^\infty\big((0,\infty)\times\X\big)$ such that
    \begin{align*}
      &\lim_{n\to\infty}\phi_n(t,x) = \mathbbm{1}_{[t_1,t_2]}(t)\phi(x) &&\text{pointwise},\\
      &\lim_{n\to\infty}\dt\phi_n(t,x) = \big(\delta_{t_1}(t)-\delta_{t_2}(t)\big)\phi(x) &&\text{weakly as measures}
    \end{align*}
    for given $t_2> t_1>0$ and a function $\phi\in C_c(\X)$.
     Because $u\in C\big([0,\infty),L^1_{\loc}(\X)\big)$, the function $t\mapsto \int_\X u(t,x)\phi(x)\,dx$ is continuous.
   Then, due to the very weak formulation~\eqref{pns:weak} and the fact that $u\in BV(\X)$, we have
   \begin{multline}\label{eq:l1-test}
    0 = \lim_{n\to\infty}\int_0^\infty\int_\X u(t,x)\dt \phi_n(t,x)-\big(\Lap_uu\big)(t,x)\phi_n(t,x)\,dx\,dt\\
    = \int_\X \big(u(t_1,x)-u(t_2,x)\big)\phi(x)\,dx - \int_{t_1}^{t_2}\int_\X\big(\Lap_uu\big)(t,x)\phi(x)\,dx\,dt.
   \end{multline}
   The assumed continuity also allows us to approach the case $t_1 = 0$.
   
   Because $u\in L^\infty\big([0,\infty),BV(\X)\big)$,
   for every $t\geq 0$ we have $u(t)\in L^1(\X)$ and $\Lap_{u(t)}u(t)\in L^1(\X)$.
   Since equality~\eqref{eq:l1-test} holds for every function $\phi\in C_c(\X)$,
   it follows that
   \begin{equation*}
    u(t_2,x) = u(t_1,x) + \int_{t_1}^{t_2}-\big(\Lap_uu\big)(t,x)\,dt\quad\text{a.e. in $x$}.
   \end{equation*}
   By~\cite[Theorem 1.4.35]{MR1691574} we obtain $u\in W^{1,1}_{\loc}\big([0,\infty),L^1(\X)\big)$.
   \end{proof}
  \end{theorem}
  \begin{remark}\label{rem:ac}
   As stated in \cite[Theorem 1.4.35]{MR1691574}, if $u\in W^{1,1}_{\loc}\big([0,\infty),L^1(\X)\big)$
   then in particular $u:[0,\infty)\to L^1(\X)$ is absolutely continuous.
  \end{remark}

  The next theorem provides the proof of the $L^1$-contraction of strong solutions, which then directly implies uniqueness.
  \begin{theorem}\label{thm:contraction}
   If $u$ and $v$ are strong solutions to problem~\eqref{pns}, corresponding to initial conditions 
  $u_0,v_0\in  BV(\X)\cap L^\infty(\X)$ respectively, then
   \begin{equation*}
    \|u(t)-v(t)\|_1\leq \|u_0-v_0\|_1 \qquad \text{for every} \quad t>0.
   \end{equation*}
   \begin{proof}
    Consider a sequence of test functions $\{\psi_n\}\subset C_c^\infty\big([0,\infty)\times\X\big)$  such that
    $\lim_{n\to\infty}\psi_n(t,x) = \mathbbm{1}_{[t_1,t_2]}(t)\sgn\big(u(t,x)-v(t,x)\big)$ for almost every $t > 0$ and $x\in\X$ and some $t_2>t_1>0$.
    Then for every $t\in[t_1,t_2]$ we have
    \begin{multline*}
    \lim_{n\to\infty}-\int_\X \big (u(t,x)-v(t,x)\big)\dt \psi_n(t,x)\,dx \\
    = \lim_{n\to\infty}\int_\X \dt \big(u(t,x)-v(t,x)\big) \psi_n(t,x)\,dx = \int_\X \dt \big|u(t,x)-v(t,x)\big|\,dx
    \end{multline*}
    and 
    \begin{multline*}
    \lim_{n\to\infty}\int_\X \Big(\big(\Lap_uu\big)(t,x)-\big(\Lap_vv\big)(t,x)\Big)\psi_n(t,x)\,dx\,dt.\\
    \int_\X \Big(\big(\Lap_uu\big)(t,x)-\big(\Lap_vv\big)(t,x)\Big)\sgn\big(u(t,x)-v(t,x)\big)\,dx\,dt.
    \end{multline*}
    Then
    \begin{multline*}
    \int_{t_1}^{t_2}\int_\X \dt \big|u(t,x)-v(t,x)\big|\,dx = 
    \int_{t_1}^{t_2}\dt \int_\X \big|u(t,x)-v(t,x)\big|\,dx\,dt \\= 
    \big\|u(t_2)-v(t_2)\big\|_1-\big\|u(t_1)-v(t_1)\big\|_1.
    \end{multline*}
    Moreover, thanks to Theorem~\ref{thm:regularity} and Remark~\ref{rem:ac}, we may also approach the case
    $t_1 = 0$.
    Finally we obtain
    \begin{multline*}
    \big\|u(t_1)-v(t_1)\big\|_1-\big\|u(t_2)-v(t_2)\big\|_1
    \\=
    \int_{t_1}^{t_2}\int_\X \Big(\big(\Lap_uu\big)(t,x)-\big(\Lap_vv\big)(t,x)\Big)\sgn\big(u(t,x)-v(t,x)\big)\,dx\,dt.
    \end{multline*}
    The right-hand side is positive due to Lemma~\ref{lemma:bv} and Theorem~\ref{thm:kato}.
   \end{proof}
  \end{theorem}

\subsection*{Properties of strong solutions}

We conclude our reasoning by gathering some of the most fundamental properties
of strong solutions to problem~\eqref{pns}.

\begin{corollary}[Uniqueness of solutions]\label{cor:uniqueness}
 Let $\meas$ be a homogeneous jump kernel. For every $u_0\in BV(\X)\cap L^\infty(\X)$ there exists a unique
 strong solution to problem~\eqref{pns}.
 \begin{proof}
  It follows directly from Theorems~\ref{thm:existence} an~\ref{thm:contraction}.
 \end{proof}
\end{corollary}

 \begin{corollary}[$L^p$-estimates]\label{cor:Lp-weak}
   If $u$ is a strong solution to problem~\eqref{pns} with initial condition $u_0\in BV(\X)\cap L^\infty(\X)$ then 
   \begin{equation*}
    \|u(t)\|_p\leq\|u_0\|_p\quad\text{and}\quad \|u(t)\|_{BV}\leq\|u_0\|_{BV} \qquad  \text{for all} \quad  t>0.    
   \end{equation*}
  \begin{proof}
   We established these estimates in relations~\eqref{p:est} and \eqref{BV:est}, in the course of proving Theorem~\ref{thm:existence} for
   the solution obtained as the limit of a subsequence of approximate solutions.
   It follows from Corollary~\ref{cor:uniqueness} that there are no other strong solutions.
  \end{proof}
 \end{corollary}
 
 \begin{corollary}[Mass conservation]\label{cor:mass_conservation}
   If $u$ is a strong solution to problem~\eqref{pns} with initial condition $u_0\in BV(\X)\cap L^\infty(\X)$ then $\int_\X u(t)\,dx = \int_\X u_0\,dx$ for every $t\geq0$.
  \begin{proof}
  The technique of this proof is essentially the same as the one used in the proof of Theorem~\ref{thm:regularity}.
  Consider a sequence $\{\phi_n\}\in C_c^\infty\big((0,\infty)\times\X\big)$ such that
    \begin{align*}
      &\lim_{n\to\infty}\phi_n(t,x) = \mathbbm{1}_{[t_1,t_2]}(t)\mathbbm{1}_K(x) &&\text{pointwise},\\
      &\lim_{n\to\infty}\dt\phi_n(t,x) = \big(\delta_{t_1}(t)-\delta_{t_2}(t)\big)\mathbbm{1}_K(x) &&\text{weakly as measures}
    \end{align*}
    for given $t_2> t_1>0$.
   Then, due to the very weak formulation~\eqref{pns:weak} and the fact that $u\in BV(\X)$, we have
   \begin{multline}\label{eq:mass}
    0 = \lim_{n\to\infty}\int_0^\infty\int_\X \dt \phi_n(t,x)u(t,x)-\big(\Lap_uu\big)(t,x)\phi_n(t,x)\,dx\,dt\\
    = \int_\X \big(u(t_1,x)-u(t_2,x)\big)\,dx - \int_{t_1}^{t_2}\int_\X\big(\Lap_uu\big)(t,x)\,dx\,dt.
   \end{multline}
   The assumed continuity also allows us to approach the case $t_1 = 0$.
   
   We may now use the symmetrization argument (see Remark~\ref{rem:sym}) and get
   \begin{equation*}
    \int_\X\Lap_uu\,dx = \iint_{\XX} \big[u(t,y)-u(t,x)\big]\kernelt{u}\,dy\,dx = 0.
   \end{equation*}
   Note that the integrals are convergent because of Lemma~\ref{lemma:bv}
   and we may use the Fubini-Tonelli theorem needed for the argument to work.
   In consequence, it follows from~\eqref{eq:mass} that
   \begin{equation*}\int_\X u(t_1,x)\,dx = \int_\X u(t_2,x)\,dx.\end{equation*}
   This means that the function $t\mapsto \int_\X u(t,x)\,dx$ is constant.
  \end{proof}
 \end{corollary}
 
 \begin{corollary}[Comparison principle]\label{cor:comparison}
 If $u$ and $v$ are strong solutions to problem~\eqref{pns} with initial conditions 
 $u_0,v_0\in BV(\X)\cap L^\infty(\X)$,
 respectively, such that $u_0(x)\leq v_0(x)$ almost everywhere in $x\in\X$ then $u(t,x)\leq v(t,x)$ almost everywhere in $(t,x)\in[0,\infty)\times \X$.

  \begin{proof}
   Using the $L^1$-contraction property established in Theorem~\ref{thm:contraction} 
   and the conservation of mass from Corollary~\ref{cor:mass_conservation} we get
   \begin{multline*}
    \int_\X \big(u(t)-v(t)\big)^+\,dx = \int_\X \frac{\big|u(t)-v(t)\big|+u(t)-v(t)}{2}\,dx 
   \\ 
    \leq
   \int_\X \frac{|u_0-v_0|+u_0-v_0}{2}\,dx = \int_\X (u_0-v_0)^+\,dx.
   \end{multline*}
   In particular,  the equality $(u_0-v_0)^+=0$ a.e.~implies 
   $(u-v)^+=0$ a.e., which means exactly that $u\leq v$ a.e.~in $[0,\infty)\times\X$.
   \end{proof}
 \end{corollary}
 \begin{corollary}[Positivity of solutions]\label{cor:positivity}
  If $u$ is a strong solution to problem~\eqref{pns} with initial condition
  $u_0\in BV(\X)\cap L^\infty(\X)$
  such that $u_0\geq 0$ then $u\geq 0$.
  \begin{proof}
   This follows directly from Corollary~\ref{cor:comparison} and the fact that $v\equiv0$ is a strong solution
   to problem~\eqref{pns} (cf. Remark~\ref{rem:constant}).
  \end{proof}
 \end{corollary}
 \subsection*{Existence of solutions for less regular initial data}
 The results obtained for strong solutions may be used to show existence of solutions
 for more general initial conditions, namely in the space $\Lpp\X$ (see Corollary~\ref{cor:main}).
 \begin{theorem}\label{thm:existence2}
  If $\meas$ is a homogeneous jump kernel then there exists a very weak solution to problem~\eqref{pns} for every
  initial condition $u_0\in \Lpp\X$.
  \begin{proof}
   We define $u_0^\epsilon = \omega_\epsilon * u_0$ for a sequence of mollifiers $\{\omega_\epsilon\}$.
   We then have $u_0^\epsilon\in BV(\X)$,
   $\|u_0^\epsilon\|_\infty\leq \|u_0\|_\infty$ and
   $\lim_{\epsilon\to0}\|u_0^\epsilon-u_0\|_1=0$.
   
   For every $\epsilon>0$ we consider the strong solution $u^\epsilon$ corresponding to
   the initial condition $u_0^\epsilon$.
   It follows from Theorem~\ref{thm:contraction} that 
   \begin{equation*}
    \sup_{t\geq0}\|u^{\epsilon_1}(t)-u^{\epsilon_2}(t)\|_1\leq \|u_0^{\epsilon_1}-u_0^{\epsilon_2}\|_1,
   \end{equation*}
   which shows that $\{u^\epsilon\}$ is a Cauchy sequence in the space
   $C_b\big([0,\infty),L^1(\X)\big)$ and hence has a limit $u\in C_b\big([0,\infty),L^1(\X)\big)$.
   
   Consider the following equations (cf. Definition~\ref{def:distributional})
   \begin{multline}\label{eq:approx2}
    \int_0^\infty\int_\X u^\epsilon(t,x)\dt \psi(t,x)\,dx\,dt 
    = \int_0^\infty\int_\X u(t,x) \big(\Lap_{u^\epsilon}\psi\big)(t,x)\,dx\,dt \\
    = \int_0^\infty\iint_{\XX} u^\epsilon(t,x)\big[\psi(t,x)-\psi(t,y)\big]\meas_{u^\epsilon(t),x,y}\,dy\,dx\,dt = 0.
   \end{multline}
   Notice that because $u$ is the limit of $\{u^\epsilon\}$ in the space $C_b\big([0,\infty),L^1(\X)\big)$,
   it also is its limit in the space $L^1\big([0,T]\times\X\big)$ for an arbitrary $T>0$.
   It follows from \cite[Theorem~4.9]{MR2759829} that there exists a function $v\in L^1\big([0,T]\times\X\big)$ and a subsequence $\epsilon_j$ such that
   \begin{align*}
   &\lim_{j\to\infty} u^{\epsilon_j}(t,x) = u(t,x)& &\text{a.e.~in $(t,x)\in[0,T]\times\X$,}\\
   &|u^{\epsilon_j}(t,x)|\leq v(t,x)& &\text{a.e.~in $(t,x)\in[0,T]\times\X$ for every $j\in\mathbb{N}$}.
   \end{align*}
   We thus have
   \begin{multline*}
    \big|u^\epsilon(t,x)\big[\psi(t,x)-\psi(t,x-y)\big]\meas_{u^\epsilon(t),x,x-y}\big|\\
    \leq 2 v(t,x)\sup_{t\in[0,T]}\|\psi(t)\|_{W^{1,\infty}(\X)}\big(1\wedge|y|\big)m_R\big(|y|\big),
   \end{multline*}
   where $R\geq \|u_0\|_\infty$. This allows us to pass to the limits on both sides of equations~\eqref{eq:approx2}
   on a subsequence $\{\epsilon_j\}$
   and verify that $u$ satisfies Definition~\ref{def:distributional}.
  \end{proof}
 \end{theorem}
 \begin{corollary}\label{cor:properties}
  Let $u$ be the very weak solution to problem~\eqref{pns} constructed in Theorem~\ref{thm:existence2}.
  We have
  \begin{itemize}
  \item $u\in C\big([0,\infty),L^1(\X)\big)$
  \item $\int_{\X} u(t,x)\,dx = \int_{\X} u_0(x)\,dx$ for all $t\geq 0$;
  \item $\|u(t)\|_p\leq \|u_0\|_p$ for all $p\in [1,\infty]$ and $t\geq 0$;
  \item if $u_0(x)\geq 0$ for almost every $x\in\X$ then $u(t,x)\geq 0$ for almost every $x\in\X$ and $t>0$.
  \end{itemize}
  \begin{proof}
   The first claim follows from the construction itself. Two other are a consequence of the fact
   that $u$ is the pointwise limit of the sequence of approximations, for which these
   claims are satisfied.
  \end{proof}
 \end{corollary}

 \end{section}
 
\begin{section}{Examples}\label{sec:examples}
In this section we discuss several examples of homogeneous jump kernels, either well-known or new.
Before we do so, however, we would like to have a quick look at the geometry of the set of jump kernels.

\begin{lemma}\label{lemma:cone}
 Let $\meas_1$ and $\meas_2$ be homogeneous jump kernels in the sense of Definition~\ref{def:jump-kernel}. Then $\meas=\alpha\meas_1+\beta\meas_2$ is a homogeneous jump kernel
 for every $\alpha,\beta\geq0$ (i.e.~the set of homogeneous jump kernels is a convex cone).
 \begin{proof}
  It is easy to see that $\meas$ satisfies conditions~\ref{a1:positive}, \ref{a2:symmetry}, \ref{a4:homogeneity} and~\ref{a6:lipschitz}. Then,
  \begin{multline*}
   (a-b)\kernell{a}{b} = \alpha(a-b)\meas_1(a,b;x,y)+\beta(a-b)\meas_2(a,b;x,y)\\
   \geq \alpha(c-d)\meas_1(c,d;x,y)+\beta(c-d)\meas_2(c,d;x,y) = (c-d)\kernell{c}{d},
  \end{multline*}
  which confirms~\ref{a3:monotonicity}, and
  \begin{multline*}
   \sup_{-R\leq a,b\leq R}\kernell{a}{b}\leq 
   \alpha \sup_{-R\leq a,b\leq R} \meas_1(a,b;x,y) + \beta\sup_{-R\leq a,b\leq R}\meas_2(a,b;x,y)\\
   \leq \alpha\, m^1_R\big(|x-y|\big)+\beta\, m^2_R\big(|x-y|\big) = m_R\big(|x-y|\big),
  \end{multline*}
  where $m_R^1$ and $m_R^2$ are functions related to $\meas_1$ and $\meas_2$ as in condition~\ref{a5:1levy}, respectively, and it follows
  that
  \begin{equation*}
   \int_\X \big(1\wedge|y|\big)m_R\big(|y|\big)\,dy 
   \leq \alpha K^1_R + \beta K^2_R,
  \end{equation*}
  where $K_R^1$ and $K_R^2$ are appropriate constants related to $m_R^1$ and $m_R^2$ as in condition~\ref{a5:1levy}.
  This confirms condition~\ref{a5:1levy} in case of the jump kernel $\meas$.
 \end{proof}
\end{lemma}
\begin{remark}\label{rem:cutout}
 It is easy to observe that if $\meas$ is a homogeneous jump kernel and $m\in L^\infty([0,\infty))$ satisfies $m\geq 0$, then $m\big(|x-y|\big)\meas\big(a,b,|x-y|\big)$ is a ho\-mo\-ge\-neous
 jump kernel as well. As a useful example we may consider $m(z)=\mathbbm{1}_{z<\delta}(z)$.
\end{remark}

\subsection*{Decoupled jump kernels}
  In this subsection we study examples of homogeneous jump kernels given in the following form
  \begin{equation}\label{separated}
   \kernell{a}{b} = F(a,b)\times \mu\big(|x-y|\big),
  \end{equation}
  where $F\geq0$ and $\mu$ is a density of a Lévy measure with low singularity, i.e.
  \begin{equation*}
   \int_\X\big(1\wedge|y|\big)\mu\big(|y|\big)\,dy\leq \infty.
  \end{equation*}
  We call such jump kernels \emph{decoupled}. We are going to consider
  several different choices of $F$ and show that the resulting functions possess
  properties~\ref{a1:positive} to~\ref{a6:lipschitz}, confirming they are homogeneous jump kernels.
  Some of these examples have been well-studied before, but we are able to verify that they fit
  neatly in our framework.
  Because of the structure of formula~\eqref{separated}, condition~\ref{a4:homogeneity} is always
  satisfied. The same is true for~\ref{a5:1levy}, as long as
  $\sup_{-R\leq a,b\leq R}F(a,b)<\infty$. In identical fashion as in Lemma~\ref{lemma:cone}, we may see that decoupled jump kernels also form a convex cone on their own.
  \subsubsection*{Fractional porous medium equation}
  In our first example we show that the theory we developed may be applied to equation~\eqref{fpme},
  for restricted ranges of parameters $s$ and $m$, and some of its generalisations.
  \begin{proposition}\label{lemma:fpme}
   Let $f\in C^1(\R)$ be a non-decreasing function.
   If
   \begin{equation*}
    F(a,b) = \frac{f(a)-f(b)}{a-b},\qquad F(a,a) = f'(a)
  \end{equation*}
  then $\meas$ given by formula~\eqref{separated} is a homogeneous jump kernel.
   \begin{proof}
    Because $f\in C^1(\R)$, the function $F$ is continuous. Because the function $f$ is non-decreasing we have $F(a,a) = f'(a) \geq 0$ and 
    \begin{equation*}
    \sgn\big(f(a)-f(b)\big)= \sgn(a-b)\quad\text{or}\quad f(a)=f(b)\quad \text{for all $a\neq b$}, 
    \end{equation*}
    therefore both~\ref{a1:positive} and~\ref{a2:symmetry} are satisfied.
    Because of the assumed monotonicity of $f$, for $a\geq c\geq d\geq b$ we have $f(a)\geq f(c)\geq f(d)\geq f(b)$, hence
    \begin{equation*}
     (a-b)\frac{f(a)-f(b)}{a-b} = f(a)-f(b)\geq f(c)-f(d) = (c-d)\frac{f(c)-f(d)}{c-d}
    \end{equation*}
    and~\ref{a3:monotonicity} is also fulfilled.

    Next, let us take $-R<a,b,c< R$ such that $|c-b|,|a-b|>\epsilon$. Then
    \begin{equation*}
     \bigg|\frac{f(a)-f(b)}{a-b}-\frac{f(c)-f(b)}{c-b}\bigg| 
     \leq \frac{2}{\epsilon^{2}}\Big(\max_{|\xi|<R}|f(\xi)|+R\max_{|\xi|<R}|f'(\xi)|\Big)|a-c|.
    \end{equation*}
    This proves the local Lipschitz-continuity condition~\ref{a6:lipschitz} (see Remark~\ref{rem:lip-sym}).
   \end{proof}
  \end{proposition}
  \begin{remark}
   This example, for $\mu\big(|y|\big) = |y|^{-N-\alpha}$, $\alpha\in(0,1)$ and $f(u)=u|u|^{m-1}$, corresponds
   to the following operator
   \begin{equation*}
    \Lap_uu = \Delta^{\alpha/2}\big(u|u|^{m-1}(x)\big),
   \end{equation*}
   introduced in Section~\ref{sec:introduction} by formula~\eqref{fpme}.
   It has been thoroughly studied in~\cite{MR2954615} (for the full range $\alpha\in(0,2)$ and $m>0$).
   Even more general non-linear operators of this type, involving linear operators represented by formula~\eqref{L:lin}, were considered in~\cite{MR3724879, MR3570132}.
  \end{remark}
  \subsubsection*{Convex diffusion operator}
  In the next example we introduce a decoupled jump kernel which has not been previously studied.
\begin{proposition}\label{lemma:breslau}
 Let $f:\R\to[0,\infty)$ be a convex function and
 \begin{equation*}
  F(a,b) = f(a)+f(b).
 \end{equation*}
 Then $\meas$ given by formula~\eqref{separated} is a homogeneous jump kernel.
 \begin{proof}
  Conditions~\ref{a1:positive} and~\ref{a2:symmetry} are easy to verify.
  
  Let $a\geq c\geq d\geq b$ and take $t\in[0,1]$ such that $d=ta+(1-t)b$. Then
  \begin{multline*}
   f(a)+f(b) \geq (1-t^2)f(a)+(1-t)^2f(b)\\
   =(1-t)\big((1+t)f(a)+(1-t)f(b)\big)
   =(1-t)\big(f(a)+tf(a)+(1-t)f(b)\big).
  \end{multline*}
  Because we assume $f$ to be convex, we have
  \begin{equation*}
   tf(a)+(1-t)f(b) \geq f\big(ta+(1-t)b\big)= f(d),
  \end{equation*}
  hence
  \begin{equation*}
   f(a)+f(b)\geq (1-t)\big(f(a)+f(d)\big).
  \end{equation*}
  In the same fashion, by taking $s\in[0,1]$ such that $c=(1-s)a+sd$, we can show that
  \begin{equation*}
   f(a)+f(d)\geq (1-s)\big(f(c)+f(d)\big).
  \end{equation*}
  We also have $(1-t)(a-b) = a-\big(ta+(1-t)b\big) = a-d$ and $(1-s)(a-d) = (c-d)$, therefore
  \begin{equation*}
   (a-b)\big(f(a)+f(b)\big) \geq (1-s)(1-t)(a-b)\big(f(c)+f(d)\big)= (c-d)\big(f(c)+f(d)\big).
  \end{equation*}
  This verifies the property~\ref{a3:monotonicity}. Since every convex function is locally Lipschitz, condition~\ref{a6:lipschitz} is satisfied and hence $\meas$ is a homogeneous jump kernel.
 \end{proof}
\end{proposition}
\begin{remark}
 This example is somewhat similar to the fractional porous medium case. Indeed, consider
 $g(a) = a|a|^{m}$ for $m=2k$ and $k\in\mathbb{N}$. Then 
 \begin{equation*}
  g(a)-g(b) = (a-b)\big(|a|^{m}+|a|^{m-1}|b|+\ldots+|a||b|^{m-1}+|b|^{m}\big)  
 \end{equation*}
 and one may be tempted to ``abbreviate'' the expression on the right-hand side to obtain
 \begin{equation}\label{eq:breslau}
  G(a,b) = (a-b)(|a|^{m}+|b|^{m}).
 \end{equation}
 However, the operator $\Lap$ based on the second kernel cannot be expressed as a linear operator acting on a non-linear
 transformation of the solution, as in the case of the fractional porous medium equation.
 
 Notice that the jump kernel related to expression~\eqref{eq:breslau} satisfies the hypothesis of Proposition~\ref{lemma:breslau} for every 
 real number $m\geq 1$.
\end{remark}
\subsubsection*{Fractional $p$-Laplacian}
Our results may also be applied to equation~\eqref{plap}.
\begin{proposition}\label{lemma:fpl}
Let $\Phi:\R\to\R$ be a continuous, non-decreasing function satisfying
   $\Phi(z)\geq 0$ for $z\geq 0$, $\Phi(-z)=-\Phi(z)$ such that ${\Phi(z)}$ is locally Lipschitz-continuous
   and $\lim_{z\to0}\tfrac{\Phi(z)}{z}<\infty$. If
  \begin{equation*}
    F(a,b) = \frac{\Phi(a-b)}{a-b}
  \end{equation*}
  then $\meas$ given by formula~\eqref{separated} is a homogeneous jump kernel.
   \begin{proof}
    Because $\sgn\big(\Phi(z)\big) = \sgn(z)$ and $\Phi(-z)=-\Phi(z)$, conditions~\ref{a1:positive} and~\ref{a2:symmetry}
    are satisfied.
    We assume the limit $\lim_{z\to0}\tfrac{\Phi(z)}{z}$ to exist, which verifies condition~\ref{a5:1levy}.
    Then, for $a\geq c\geq d\geq b$ and by using the fact that 
    $\Phi$ is non-decreasing, we have
    \begin{equation*}
     (a-b)\frac{\Phi(a-b)}{a-b} = \Phi(a-b) \geq \Phi(c-d) = (c-d)\frac{\Phi(c-d)}{c-d},
    \end{equation*}
    from which~\ref{a3:monotonicity} follows. Let $-R<a,b,c<R$ such that $|a-b|,|c-b|>\epsilon$. The
    function $\Phi$ is locally Lipschitz-continuous, therefore
    \begin{equation*}
     \bigg|\frac{\Phi(a-b)}{a-b}-\frac{\Phi(c-b)}{c-b}\bigg|\leq \frac{c_R}{\epsilon^2}|a-c|,
    \end{equation*}
    with a number $c_R>0$ depending on $R$ and the Lipschitz constant of the function~$\Phi$.
    This confirms condition~\ref{a6:lipschitz} (see Remark~\ref{rem:lip-sym}).
    \end{proof}
  \end{proposition}
  \begin{remark}
   Consider the function $\Phi(z) = |z|^{p-2}z$. It is easy to verify that it satisfies the hypothesis of Proposition~\ref{lemma:fpl}
   if and only if $p\geq 2$. Then we may take $s<\frac{1}{p}$ (so that $sp<1$
   and condition~\ref{a5:1levy} holds) and we recover the non-local
   $s$-fractional $p$-Laplace operator which appears in equation~\eqref{plap}.
  \end{remark}
  \subsubsection*{Doubly non-linear Lévy operator}
  The following non-local counterpart of equation~\eqref{eq:doubly-nonlinear}, which appears to have not
  yet been studied, turns out to be covered by our theory.
  \begin{proposition}\label{lemma:doubly-nonlinear}
   Suppose functions $f$ and $\Phi$ satisfy adequate parts of hypotheses of Propositions~\ref{lemma:fpme} and~\ref{lemma:fpl}, respectively.
   If
   \begin{equation*}
    F(a,b) = \frac{\Phi\big(f(a)-f(b)\big)}{a-b}
  \end{equation*}
  then $\meas$ given by formula~\eqref{separated} is a homogeneous jump kernel.
   \begin{proof}
    Conditions~\ref{a1:positive}, \ref{a2:symmetry}, \ref{a3:monotonicity} and~\ref{a6:lipschitz} follow easily
    as combinations of the arguments already used to prove Propositions~\ref{lemma:fpme} and~\ref{lemma:fpl}.
   \end{proof}
  \end{proposition}
\subsection*{Entangled jump kernels.}
Homogeneous jump kernels which cannot be decomposed as in formula~\eqref{separated} are also part of
our framework.
\subsubsection*{Jump kernel with variable order.}
In the last example we study an operator whose ``differentialbility order'' is allowed to depend on the solution itself.

 \begin{proposition}\label{lemma:bielefeld}
  Let 
  \begin{equation*}
   \Psi(\mathfrak{a};z) = \Psi_1(\mathfrak{a})\mathbbm{1}_{z<1}(z)+\Psi_2(\mathfrak{a})\mathbbm{1}_{z\geq 1}(z)+\Theta\big(z\big),
  \end{equation*}
  where
  $\Psi_1:[0,\infty)\to \R$ is non-decreasing,
  $\Psi_2:[0,\infty)\to \R$ is non-increasing,
  both $\Psi_1$ and $\Psi_2$ are locally Lipschitz-continuous,
  $\Theta:[0,\infty)\to\R$ is measurable and
  \begin{equation*}
   0<A_1\leq \Psi(\mathfrak{a};z)\leq A_2<1,\qquad A_1\leq \Theta(z)\leq A_2
  \end{equation*}
  Then
  \begin{equation*}
    \kernell{a}{b} = |x-y|^{-N-\Psi(|a-b|;|x-y|)}
  \end{equation*}
  is a homogeneous jump kernel.
  \begin{proof}
   Conditions~\ref{a1:positive},~\ref{a2:symmetry} and~\ref{a4:homogeneity} are easy to check. In order to
   verify~\ref{a3:monotonicity} we notice that for $\mathfrak{a}\geq \mathfrak{b}$ and $z\geq0$ we have
   \begin{equation*}
    \mathfrak{a}z^{-\Psi_1(\mathfrak{a})\mathbbm{1}_{z<1}(z)}\geq \mathfrak{b}z^{-\Psi_1(\mathfrak{b})\mathbbm{1}_{z<1}(z)},
   \end{equation*}
   because $\Psi_1$ is non-decreasing and
   \begin{equation*}
    \mathfrak{a}z^{-\Psi_2(\mathfrak{a})\mathbbm{1}_{z\geq1}(z)}\geq \mathfrak{b}z^{-\Psi_2(\mathfrak{b})\mathbbm{1}_{z\geq1}(z)},
   \end{equation*}
   because $\Psi_2$ is non-increasing. Hence for $a\geq c\geq d\geq b$, by setting $\mathfrak{a} = a-b$ and $\mathfrak{b} = c-d$, we obtain
   \begin{equation*}
   (a-b)\kernell{a}{b}\geq (c-d)\kernell{c}{d}.
   \end{equation*}
   
   Notice the following estimate
   \begin{multline}\label{eq:bielefeld1}
    \int_\X\big(1\wedge|x-y|\big)\sup_{-R<a,b<R}\kernell{a}{b}\,dy \\\leq
    \int_{|y|<1}|y|^{1-N-A_2}\,dy +
    \int_{|y|\geq 1}|y|^{-N-A_1}\,dy < \infty
   \end{multline}
   and let
   \begin{equation*}
   m^0(z) = \mathbbm{1}_{z<1}(z)z^{-N-A_2}+\mathbbm{1}_{z\geq1}(z)z^{-N-A_1}.
   \end{equation*}
   Because functions $\Psi_1$, $\Psi_2$ are locally Lipschitz-continuous, we have
   \begin{equation}\label{eq:bielefeld2}
   \begin{split}
    \big|&z^{-N-\Psi(\mathfrak{a};z)}-z^{-N-\Psi(\mathfrak{b};z)}\big|\\
    &= \big(z^{-N-\Theta(z)}\big)\Big(\big|z^{\Psi_1(\mathfrak{a})}-z^{\Psi_1(\mathfrak{b})}\big|\mathbbm{1}_{z<1}(z)
    + \big|z^{\Psi_2(\mathfrak{a})}-z^{\Psi_2(\mathfrak{b})}\big|\mathbbm{1}_{z\geq 1}(z)\Big)\\
    &\leq \big|\Psi_1(\mathfrak{a})-\Psi_1(\mathfrak{b})\big|\big|\log(z)\big|
      \Big(\sup_{\subalign{\alpha&\in[0,A_2-A_1]\\z&\in[0,1)}}z^\alpha\Big)\big(z^{-N-A_2}\mathbbm{1}_{z<1}(z)\big)\\
    &\qquad+ \big|\Psi_2(\mathfrak{a})-\Psi_2(\mathfrak{b})\big|\big|\log(z)\big|
      \Big(\sup_{\subalign{\alpha&\in[A_1-A_2,0]\\z&\in[1,\infty)}}z^\alpha\Big)\big(z^{-N-A_1}\mathbbm{1}_{z\geq 1}(z)\big)\\
    &\leq\big(L_{2R}^{\Psi_1}+L_{2R}^{\Psi_2}\big)|\mathfrak{a}-\mathfrak{b}|\big|\log(z)\big|m^0(z).
     \end{split}
   \end{equation}
   In estimate~\eqref{eq:bielefeld2} we used the local Lipschitz-continuity of functions $\Psi_1$ and $\Psi_2$ as well as the fact that
   \begin{align*}
    &\sup_{\subalign{\alpha&\in[0,A_2-A_1]\\z&\in[0,1)}}z^\alpha =\sup_{\subalign{\alpha&\in[A_1-A_2,0]\\z&\in[1,\infty)}}z^\alpha = 1.
   \end{align*}
   It follows from estimates~\eqref{eq:bielefeld1} and~\eqref{eq:bielefeld2} that by putting $\mathfrak{a} = |a-b|$, $\mathfrak{b} = |c-d|$ and $z = |x-y|$ we may verify
   both conditions~\ref{a5:1levy} and~\ref{a6:lipschitz}
   by considering
   \begin{equation*}
    m_R\big(z\big) = \Big(\mathbbm{1}_{z<1}(z)z^{-N-A_2}+\mathbbm{1}_{z\geq1}(z)z^{-N-A_1}\Big)
    \times\max\big\{1,\big|\log(z)\big|\big\}.\qedhere
   \end{equation*}
  \end{proof}
  \end{proposition}
  \begin{remark}
   In Proposition~\ref{lemma:bielefeld} it would suffice to assume local Lipschitz-con\-ti\-nuity of functions $\Psi_1$ and $\Psi_2$ on $(0,\infty)$.
  \end{remark}

\end{section}
\bibliographystyle{siam}
\bibliography{kkk-06}
\end{document}